\documentclass[10pt]{article}
\usepackage{latexsym}
\usepackage{amssymb}
\usepackage{amscd}
\usepackage{amsmath}
\usepackage{array}
\usepackage{algorithm}
\usepackage{algorithmic}

\usepackage[usenames,dvipsnames]{color}
\usepackage{graphicx}
\usepackage{wrapfig}
\usepackage{tikz}
\usetikzlibrary{shapes.geometric}
\usepackage{mathtools}

\usepackage{hyperref}
\hypersetup{
 colorlinks,
 citecolor=Violet,
 linkcolor=Blue,
 urlcolor=Blue}

\usepackage{cleveref}

\newcommand{\goto}{\rightarrow}

\newcommand{\bR}{\mbox{$\mathbb{R}$}}
\newcommand{\bS}{\mbox{$\mathbb{S}$}}
\newcommand{\bZ}{\mbox{$\mathbb{Z}$}}

\newcommand{\bC}{\mbox{$\mathbb{C}$}}
\newcommand{\bitem}{\begin{itemize}}
\newcommand{\eitem}{\end{itemize}}
\newcommand{\beq}{\begin{equation}}
\newcommand{\eeq}{\end{equation}}
\newcommand{\bea}{\begin{eqnarray}}
\newcommand{\eea}{\end{eqnarray}}
\newcommand{\bear}{\begin{eqnarray*}}
\newcommand{\eear}{\end{eqnarray*}}
\newcommand{\bfig}{\begin{figure}[ht]}
\newcommand{\efig}{\end{figure}}
\newcommand{\argmin}{\mathop{\rm arg min}}
\newtheorem{theorem}{Theorem}[section]
\newtheorem{lemma}[theorem]{Lemma}

\newtheorem{proposition}[theorem]{Proposition}
\newtheorem{definition}[theorem]{Definition}

\newtheorem{remark}[theorem]{Remark}

\usepackage[mathscr]{euscript}

\newcommand{\gap}{\vspace{0.1in}}

\newcommand{\GD}{\textsf{OGD}}

\usepackage{mathtools}

\def\qed{\hfill\rule{2.0mm}{2.0mm}}
\def\pf{\noindent{\bf Proof:}~ }
\def\eop{\hfill\rule{2.0mm}{2.0mm}}

\DeclareMathOperator*{\essinf}{ess\,inf}

\begin{document}
\title{
{\sc The Essential Best and Average Rate of Convergence of the Exact Line Search Gradient Descent Method}
}
\date{April 28, 2023 \\ Revised May 10, 2025}

\author{
Thomas Yu\thanks{
Department of Mathematics, Drexel University. Email: \href{mailto:pty23@drexel.edu}{pty23@drexel.edu}.
He is supported in part by the National Science Foundation grants DMS 1522337 and DMS 1913038.
}
}

\makeatletter \@addtoreset{equation}{section} \makeatother
\maketitle

\centerline{\bf Abstract:}
 It is very well known that when the exact line search gradient descent method is applied to a convex
quadratic objective, the worst-case rate of convergence (ROC), among all seed vectors,
deteriorates as the condition number of the Hessian of the objective grows.
By an elegant analysis due to H. Akaike, it is generally believed -- but not proved -- that in the ill-conditioned regime
the ROC for almost all initial vectors, and hence also the average ROC, is close to the worst case ROC. We complete
Akaike's analysis by determining the \emph{essential best case ROC} (defined in a measure-theoretic way) by using a
dynamical system approach, facilitated by the theorem of center and
stable manifolds. Our analysis also makes apparent the effect of an intermediate eigenvalue in the Hessian by establishing
the following amusing result: In the absence of an intermediate eigenvalue, the average ROC gets arbitrarily \emph{fast} -- not slow -- as the Hessian gets increasingly ill-conditioned.

We discuss in passing some contemporary applications of exact line search GD to well-conditioned polynomial optimization problems arising from imaging and data sciences.
In particular, we observe that a tailored exact line search GD algorithm for a POP arising from the phase retrieval problem is only 50\% more expensive per iteration than its constant step size counterpart, while promising a ROC only matched by the optimally tuned (constant) step size which can rarely be achieved in practice.

\vspace{.2in} \noindent {\bf Acknowledgments.} This work was inspired by a discussion on
polynomial optimization problems with the author of \cite{Beck:Book1} during his visit to Tel Aviv University in January 2020.
He also thanks Tom Duchamp for many insightful discussions, and Kyle Steppe for help on extending Algorithm~\ref{alg:ExactLineSearchPhaseRetrieval} to the setting of coded diffraction patterns.

\vspace{.2in}


\section{Introduction}
This paper revisits an elegant but largely forgotten mathematical analysis of a rarely used version of a widely used method, namely,
the exact line search gradient descent method. 

 It is well-known that GD methods 
 converge slowly for ill-conditioned problems; and many alternative methods have been proposed to accelerate the convergence in the ill-conditioned regime. For at least those trained in the computational mathematics community circa the 70's-90's, GD was usually introduced as the first numerical method for optimization, but quickly dismissed as being a slow method, and the curriculum moved on to methods with a faster rate of convergence, such as Krylov subspace or quasi-Newton methods.
The optimization toolbox of Matlab does not even have a GD solver. Unsurprisingly, it was rare to see anyone in the traditional scientific computing community who used GD in research.

In modern AI and machine learning research, on the contrary, GD is considered the bedrock algorithm,
most notably in the training of (deep) neural networks.
This work is instead motivated by imaging and data science problems such as phase retrieval and matrix completion, which give rise to
well-conditioned unconstrained optimization problems with a polynomial structure. For these problems, GD methods appear to be effective;
see Section~\ref{sec:POP} and the references therein. Moreover,
the simplicity of GD facilitates statistical analysis, leading to useful insights.

Exact line search is usually deemed impractical in the optimization literature, but when the objective function has a specific global structure then its use can be beneficial.
A notable case is when the objective is a polynomial. Polynomial optimization problems (POPs) abound in diverse
applications; see, for example, \cite{Nie2009,MR4099988,6563125,Ahmadi2016,MR1968126,MohanYipYu:Bilayer} and Section~\ref{sec:POP}.
It is well-known that POPs can be convexified into SDPs, see \cite{doi:10.1137/S1052623400366802,MR3469431,doi:10.1137/1.9781611972290}. However, the resulted SDPs are usually prohibitively expensive to solve even with state-of-the-art interior point methods. The work of Baurer-Monteiro \cite{Burer-2003,Burer-2005}, in turn, suggests to `de-convexify' SDPs back
to quadratically constrained quadratic programs and solve them using traditional gradient based methods.
These interesting developments for general POPs, together with
extensive developments for the specific POPs arising from data and imaging sciences \cite{MR3332993,Candes:PhaseLift,CANDES2015277,MR3628877,SunQuWright,MR4099988,MR1968126,NIPS2016_7fb8ceb3,CandesRecht:MatrixCompletion},
prompt us to reconsider the theory and application of the exact line search gradient descent method.

\subsection{Background materials}
\emph{Exact line search} refers to the `optimum' choice of step size, namely
$s = \argmin_t f( x + t d)$, where $d$ is the search direction, hence the old nomenclature
\emph{optimum gradient descent} in the case of $d=-\nabla f(x)$.\footnote{We use the
terms `optimum GD' and `GD with exact line search' interchangeably in this article. The former terminology is abandoned in the contemporary
optimization literature; we use it in this article only because Akaike's original article uses it.}
When $f$ is, say, a degree 4 polynomial (see Section~\ref{sec:POP}),
it amounts to determining the univariate quartic polynomial $p(t) = f( x - t d)$ followed by finding its minimizer,
which in specific important applications can be computed efficiently, potentially leading to a practical advantage of exact line search.

Let us recall the key theoretical benefit of using exact line search in GD.
We focus on the case when the objective is a strictly convex quadratic, which locally approximates any smooth objective function
in the vicinity of a local minimizer with a positive definite Hessian.
By the invariance of GD (constant step size or exact line search) under rigid transformations,
there is no loss of generality, as far as the study of ROC is concerned, to consider only quadratics of the form
\bea \label{eq:Quadratic}
f(x)=\frac{1}{2} x^T A x, \mbox{ with } A = {\rm diag}(\lambda), \;\; \lambda = (\lambda_1,\ldots,\lambda_n), \;\; \lambda_1 \geq \cdots \geq \lambda_n>0.
\eea
Its gradient is Lipschitz continuous with Lipschitz constant $L = \lambda_1$.
Also, it is strongly convex with the strong convexity parameter $\sigma = \lambda_{n}$.
In the case of constant step size GD, we have
 $x^{(k+1)} = (I - s A) x^{(k)}$, so the rate of convergence is
 given the spectral radius of $I-sA$, which equals $\max_{1\leq i \leq n} |1-s\lambda_i|$. From this, it is easy to check that the step size
\bea \label{eq:OptimalConstant}
s = 2/(\lambda_{1}+\lambda_{n})
\eea
minimizes the spectral radius of $I-sA$, and hence offers the optimal ROC
\bea \label{eq:optimalROC}
\|x^{(k)}\| = O(\rho^k) \mbox{ with } \; \rho=\frac{\lambda_{1}-\lambda_{n}}{\lambda_{1}+\lambda_{n}}.
\eea

Gradient descent with exact line search involves \emph{non-constant} step sizes:
$x^{(k+1)} = (I - s_k A) x^{(k)}$, with $s_k = (x^{(k)})^T A^2 x^{(k)}/(x^{(k)})^T A^3 x^{(k)}$.
For convenience, denote the
iteration operator by $\GD$, i.e.
\bea \label{eq:GD}
x^{(k+1)} = \GD(x^{(k)}), \quad \GD(x) := \GD(x; \lambda) := x - \frac{x^T A^2 x}{x^T A^3 x} A x.
\eea
We set $\GD(0)=0$ so that $\GD$ is a well-defined self-map on $\bR^n$.
By norm equivalence, one is free to choose any norm in the study of ROC;
and the first trick is to notice that by choosing the $A$-norm, defined by $\|x\|_A := \sqrt{x^T A x}$,
we have the following convenient relation:
\bea \label{eq:GDnorm}
\| \GD(x) \|_A^2 = \left[1-\frac{(x^T A^2 x)^2}{(x^T A x)(x^T A^3 x)} \right] \|x\|_A^2.
\eea
Write $d = A x$ ($= \nabla f(x)$). By the Kantorovich inequality,\footnote{For a proof of the Kantorovich inequality $\frac{(d^T d)^2}{(d^T A^{-1} d)(d^T A d)}
\geq \frac{4 \lambda_{1} \lambda_{n}}{(\lambda_{1} + \lambda_{n})^2} = \frac{4 {\rm cond}(A)}{(1+{\rm cond}(A))^2}$, see, for example, \cite{MR136623,Beck:Book1}. A good way to appreciate this result is to compare it with the following
bound obtained from Rayleigh quotient: $\frac{(d^T d)^2}{(d^T A^{-1} d)(d^T A d)}\geq \frac{\lambda_{n}}{\lambda_{1}} = {\rm cond}(A)^{-1}$.
Unless $\lambda_1=\lambda_n$, Kantorovich's bound is always sharper, and is the sharpest possible in
the sense that there exists a vector $d$ so that equality holds.}
\bea \label{eq:Kantorovich}
\frac{(x^T A^2 x)^2}{(x^T A x)(x^T A^3 x)} = \frac{(d^T d)^2}{(d^T A^{-1} d)(d^T A d)}
\geq \frac{4 \lambda_{1} \lambda_{n}}{(\lambda_{1} + \lambda_{n})^2},
\eea
which yields the well-known error bound for the optimum GD method:
\bea \label{eq:Worst}
\|x^{(k)}\|_A \leq \Big( \frac{\lambda_1-\lambda_n}{\lambda_1+\lambda_n} \Big)^k \|x^{(0)}\|_A.
\eea
So optimum GD satisfies the same upper-bound for its ROC as in \eqref{eq:optimalROC}.
\emph{The constant step size GD method with the optimal choice of step size \eqref{eq:OptimalConstant}, however, should not be confused
with the optimum GD method}, they have the following fundamental differences:

\begin{itemize}
\item The optimal step size \eqref{eq:OptimalConstant}
requires the knowledge of the two extremal eigenvalues, of which the determination is
no easier than the original minimization problem. In contrast, the optimum GD method is blind to the
values of $\lambda_{1}$ and $\lambda_n$.
\item Due to the linearity of the iteration process, 
GD with the optimal constant step size \eqref{eq:OptimalConstant}
achieves the ROC $\|x^{(k)}\| \asymp C \rho^k$, with $\rho$ in \eqref{eq:optimalROC},
for any initial vector with a non-zero component in the dominant eigenvector of $A$, and hence \emph{for almost all initial vectors} $x^{(0)}$.
So for this method \textbf{the worst case ROC is the same as the average ROC}. (As
the ROC is invariant under scaling of $x^{(0)}$, the average and
worst case ROC can be defined by taking the average and maximum, respectively, ROC over all seed vectors of unit length.)
In contrast, $\GD$ is nonlinear and the worst case ROC
\eqref{eq:optimalROC} is attained only for specific initial vectors $x^{(0)}$; see Proposition~\ref{prop:SharpestOneStep}. It is much less obvious how the average
ROC, defined by \eqref{eq:AveROC} below, compares to the worst case ROC.
\end{itemize}

Due to \eqref{eq:GDnorm}, we define the \textbf{(one-step) shrinking factor} by
\bea \label{eq:rho}
\rho(x,\lambda) = \sqrt{1-\frac{(x^T A^2 x)^2}{(x^T A x)(x^T A^3 x)}}
= \sqrt{1- \frac{(\sum_i \lambda_i^2 x_i^2)^2}{(\sum_i \lambda_i x_i^2)(\sum_i \lambda_i^3 x_i^2)}}.
\eea
Then, for any initial vector $x^{(0)} \neq 0$, the rate of convergence of the optimum gradient descent method
applied to the minimization of
\eqref{eq:Quadratic} is given by
\bea \label{eq:ROC_def}
\rho^\ast(x^{(0)},\lambda) := \limsup_{k\goto \infty} \Bigg[ \prod_{j=0}^{k-1} \rho(\GD^j( x^{(0)}), \lambda) \Bigg]^{1/k}.
\eea
As $\rho^\ast(x^{(0)},\lambda)$ depends only on the direction of $x^{(0)}$, and is
insensitive to sign changes in the components of $x^{(0)}$ (see \eqref{eq:invariance1}), the \textbf{average ROC}
can be defined based on averaging over all $x^{(0)}$ on the unit sphere, or just over the positive octant
of the unit sphere. Formally,
\begin{definition} The average ROC of the optimum GD method applied to \eqref{eq:Quadratic} is
\bea \label{eq:AveROC}
\mbox{Average ROC} := \int_{\bS^{n-1}} \rho^\ast(x,\lambda) d\mu(x) = 2^n \int_{\bS^{n-1}_+} \rho^\ast(x,\lambda) d\mu(x),
\eea
where $\mu$ is the uniform probability measure on $\bS^{n-1}$, and $\bS^{n-1}_+ := \{x \in \bS^{n-1}: x\geq 0\}$.
\end{definition}
We have
\bea \label{eq:AverageWorst}
\mbox{Average ROC} \leq \mbox{Worst case ROC} \stackrel{(\ast)}{=} \frac{1-a}{1+a}, \;\;\; \mbox{ where } a = \frac{\lambda_{n}}{\lambda_{1}} = {\rm cond}(A)^{-1}.
\eea
Note that \eqref{eq:Worst} only shows that the worst case ROC is upper bounded by
$(1-a)/(1+a)$; for a proof of the equality ($\ast$) in \eqref{eq:AverageWorst}, see Proposition~\ref{prop:SharpestOneStep}.

\subsection{Main results} \label{sec:Main}
In this paper, we establish the following result:
\begin{theorem} \label{thm:Main}
(i) If $A$ has only two distinct eigenvalues, then the average ROC approaches 0 
when ${\rm cond}(A) \goto \infty$.
(ii) If $A$ has an intermediate eigenvalue $\lambda_i$ uniformly bounded away from the two extremal eigenvalues, then
the average ROC approaches the worst case ROC in \eqref{eq:AverageWorst}, which approaches 1,
when ${\rm cond}(A) \goto \infty$.
\end{theorem}

We shall prove a stronger version of the second part of Theorem~\ref{thm:Main}:
\begin{theorem} \label{thm:AkaikeBound}
If $A$ has an intermediate eigenvalue, i.e. $n>2$ and there exists $i \in \{ 2,\ldots,n-1 \}$ s.t. $\lambda_{1} > \lambda_i > \lambda_{n}$,
then 
\bea \label{eq:LowerBound}
\essinf_{x^{(0)} \in \mathbb{S}^{n-1}} \rho^\ast(x^{(0)},\lambda) = \frac{1-a}{\sqrt{(1+a)^2+B a}},
\eea
where $a= {\rm cond}(A)^{-1}=\lambda_n/\lambda_1$, $B = \frac{4(1+\delta^2)}{1-\delta^2}$,
$\delta = \min_{i: \lambda_{1} > \lambda_i > \lambda_{n}} \frac{|\lambda_i-(\lambda_{1}+\lambda_{n})/2|}{(\lambda_{1}-\lambda_{n})/2}$.
\end{theorem}

Together with the standard worst-case bound \eqref{eq:Worst}, we have for almost all seed vectors $x^{(0)}$,
\beq \label{eq:LowerUpper}
\frac{1-a}{\sqrt{(1+a)^2+B a}} \leq \rho^\ast(x^{(0)},\lambda) \leq \frac{1-a}{1+a}.
\eeq
Note that when there is at least one intermediate eigenvalue $\lambda_i$ uniformly bounded away from the two extremal eigenvalues, the value $B$ above is uniformly bounded
above. (In particular, $B=4$ if there is an intermediate eigenvalue $\lambda_i$ equals to $(\lambda_{1}+\lambda_{n})/2$.)
Therefore, this result says that \emph{in the ill-conditioned regime (i.e. when $a={\rm cond}(A)^{-1}$ is small), the essential best case ROC, and hence also the average ROC, is only marginally better than the worst case ROC}. This explains why the second part of \eqref{thm:Main} follows from Theorem~\ref{thm:AkaikeBound}.


\begin{remark} \normalfont
It is shown in \cite[\S2]{Akaike-1959} that $\rho^\ast(x^{(0)},\lambda)$ is lower-bounded by the right-hand side of
\eqref{eq:LowerBound} with the proviso of a difficult-to-verify condition on $x^{(0)}$.
 The undesirable condition seems to be an artifact of
the subtle argument in \cite[\S2]{Akaike-1959}, which also makes it hard to see whether the bound \eqref{eq:LowerBound} is tight.
Our proof of Theorem~\ref{thm:AkaikeBound} uses Akaike's results in \cite[\S1]{Akaike-1959}, but replace his arguments
in \cite[\S2]{Akaike-1959} by a more natural dynamical system approach. The proof
shows that the bound is tight and holds for a set of $x^{(0)}$ of full measure, which also allows us to conclude the second part
of Theorem~\ref{thm:Main}.
It uses
the center and stable manifolds theorem, a result that was not available at the time \cite{Akaike-1959} was written.
\end{remark}

\begin{remark} \normalfont
For constant step size GD, ill-conditioning \emph{alone} is enough to cause slow convergence for almost all initial vectors.
For exact line search, however, it is ill-conditioning in cahoot with an intermediate eigenvalue that causes the slowdown. This is already apparent
from Akaike's analysis; the first part of
Theorem~\ref{thm:Main} intends to bring this point home, by showing that the exact opposite happens in the absence of an intermediate eigenvalue.
\end{remark}

\noindent
{\bf Organization of the proofs.} In Section~\ref{sec:Akaike}, we recall the results of Akaike, along with a few
observations not directly found in \cite{Akaike-1959}. These results and observations are necessary for the final proofs of Theorem~\ref{thm:Main}(i) and Theorem~\ref{thm:AkaikeBound}, to be presented in Section~\ref{sec:2D} and Section~\ref{sec:TheoremAkaikeBound}, respectively.
The key idea of Akaike is an interesting discrete dynamical system, with a probabilistic interpretation,
 underlying the optimum GD method. We give an exposition of this dynamical system in Section~\ref{sec:Akaike}. A key property of this dynamical system is recorded in Theorem~\ref{thm:AkaikeTheorem2}. In a nutshell, the theorem tells us that part of the properties of the dynamical system in the high-dimensional
 case is captured by the 2-dimensional case; so -- while the final results say that there is a drastic difference between the $n=2$ and $n>2$ cases --
 the analysis of the latter case (in Section~\ref{sec:TheoremAkaikeBound}) uses some of the results in
 the former case (in Section~\ref{sec:2D}).
Appendix~\ref{AppendixA} records two auxiliary technical lemmas. Appendix~\ref{AppendixB} recalls a version of the theorem of center and stable manifolds stated in
the monograph \cite{MR869255}, and discusses a refinement of the result needed to prove the full version of Theorem~\ref{thm:AkaikeBound}.

Section~\ref{sec:2D} and Section~\ref{sec:TheoremAkaikeBound} also present computations,
in Figures~\ref{fig:WorstvsAverage} and \ref{fig:Figure4}, that serve to
illustrate some key ideas behind the proofs.

\gap
Before proceeding to the proofs, we consider some contemporary applications of exact line search methods to
POPs. Section~\ref{sec:POP} below is independent from the rest of the paper.

\subsection{Applications of exact line search methods to POPs} \label{sec:POP}
In its abstract form, the phase retrieval problem
seeks to recover a signal $x \in \mathbb{K}^n$ ($\mathbb{K} = \bR$ or $\bC$) from its noisy `phaseless measurements'
$y_i \approx |\langle x, a_i \rangle|^2$, with enough random
`sensors' $a_i \in \mathbb{K}^n$. A plausible
approach is to choose $x$ that solves
\bea \label{eq:PhaseRetrievalPOP}
\min_{x \in \mathbb{K}^n} \sum_{j=1}^m \Big[ y_j - |\langle x, a_j \rangle|^2 \Big]^2.
\eea
The two squares makes
it a degree 4 POP.

We consider also another data science problem: matrix completion.
In this problem, we want to exploit the a priori \emph{low rank} property of a data matrix $M$
in order to estimate it from just a small fraction of its entries $M_{i,j}$, $(i,j) \in \Omega$.
If we know a priori that $M \in \bR^{m\times n}$ has rank $r \ll \min(m,n)$, then similar to \eqref{eq:PhaseRetrievalPOP} we may hope to
recover $M$ by solving
\bea \label{eq:MatrixCompletion}
\min_{X \in \bR^{m\times r}, Y \in \bR^{n \times r}} \sum_{(i,j)\in \Omega} \Big[ (XY^T)_{i,j} - M_{i,j} \Big]^2.
\eea
It is again a degree 4 POP.


Extensive theories have been developed for addressing the following questions:
(i) Under what conditions -- in particular how big the sample size $m$ for phase retrieval and $|\Omega|$ for matrix completion --
would the global minimizer of \eqref{eq:PhaseRetrievalPOP} or \eqref{eq:MatrixCompletion} recovers the underlying object of interest?
(ii) What optimization algorithms would be able to compute the global minimizer?

It is shown in \cite{MR4099988} that constant step size GD with an appropriate choice of initial vector and step size
applied to the optimization problems above probably guarantees
success in recovering the object of interest under suitable statistical models. For the phase retrieval problem,
assume for simplicity $x^\ast \in \bR^n$, $y_j=(a_j^T x^\ast)^2$,
write
\bea \label{eq:PhaseRetrievalQuartic}
f(x) := \frac{1}{4m} \sum_{j=1}^m \Big[ y_j - (a_j^T x)^2 \Big]^2.
\eea
Then
\bea \label{eq:PhaseRetrievalQuarticGrad}
\nabla f(x) = -\frac{1}{m} \sum_{j=1}^m \Big[ y_j - (a_j^T x)^2 \Big] (a_j^T x) a_j \;\mbox{ and } \;
\nabla^2 f(x) = \frac{1}{m} \sum_{j=1}^m \Big[ 3(a_j^T x)^2 -y_j \Big] a_j a_j^T.
\eea
Under the Gaussian design of sensors
$a_j \stackrel{\rm i.i.d.}{\sim} N(0,I_n)$, $1\leq j\leq m$, considered in, e.g., \cite{Candes:PhaseLift,MR3332993,MR4099988}, we have
$$
\mathbb{E} \big[ \nabla^2 f(x) \big] = 3 \big[ \|x\|_2^2 I_n + 2xx^T \big] - \big[ \|x^\ast\|_2^2 I_n + 2x^\ast (x^\ast)^T \big].
$$
At the global minimizer $x=x^\ast$, $\mathbb{E} \big[ \nabla^2 f(x^\ast) \big] = 2\big[ \|x^\ast\|_2^2 I_n + 2x^\ast (x^\ast)^T \big]$, so
$$
{\rm cond}(\mathbb{E} \big[ \nabla^2 f(x^\ast) \big]) = 3.
$$
This suggests that when the sample size $m$ is large enough,
we may expect that the Hessian of the objective is well-conditioned for $x \approx x^\ast$ -- \emph{hence the use of gradient descent is appropriate}.
Indeed, when $m \asymp n \log n$, 
the discussion in \cite[Section 2.3]{MR4099988}
implies that ${\rm cond}(\nabla^2 f(x^\ast))$ grows slowly with $n$:
$$
{\rm cond}(\nabla^2 f(x^\ast)) = O(\log n)
$$
with a high probability.
However, unlike \eqref{eq:Quadratic}, the objective \eqref{eq:PhaseRetrievalQuartic} is a quartic instead of a quadratic polynomial, so
the Hessian $\nabla^2 f(x)$ is not constant in $x$.
The problem exhibits the following phenomena:
\begin{itemize}
\item[(i)] On the one hand, in the directions given by $a_j$, the Hessian $\nabla^2 f(x^\ast + \delta a_j/\|a_j\|)$ has a condition number
that grows (up to logarithmic factors) as $O(n)$, meaning that the objective can get increasingly ill-conditioned as the dimension $n$ grows, even within a small
ball around $x^\ast$ with a fixed radius $\delta$.
\item[(ii)] On the other hand, most directions $v$ would not be
 too close to be parallel to $a_j$, and ${\rm cond} \big( \nabla^2 f(x^\ast + \delta v/\|v\|) \big) = O(\log n)$ with a high probability.
\item[(iii)] Constant step-size GD, with a step size that can be chosen nearly constant in $n$, has the property of staying away from the
ill-conditioned directions, hence no pessimistically small step sizes
or explicit regularization steps avoiding the bad directions are needed. Such an
 `implicit regularization' property of constant step size GD is the main theme of the article \cite{MR4099988}.
\end{itemize}

To illustrate (i) and (ii) numerically, we compute the condition numbers of the Hessians $\nabla^2 f(x)$ at $x=x^\ast=\bar{x}/\|\bar{x}\|$, $x=x^\ast + .5 a_1/\|a_1\|$
and $x=x^\ast + .5 z/\|z\|$ for $a_j, \bar{x}, z \stackrel{\rm i.i.d.}{\sim} N(0,I_n)$, with $n=1000 k$, , $k=1,\ldots,5$, and $m=n\log_2(n)$:
\begin{center}
\begin{tabular}{|c|c|c|c|}
  \hline
   & ${\rm cond}(\nabla^2 f(x^\ast))$ &  ${\rm cond}(\nabla^2 f(x^\ast + .5 a_1/\|a_1\|)$ &  ${\rm cond}(\nabla^2 f(x^\ast + .5 z/\|z\|))$ \\
   \hline
$n=1000$ & 14.1043 & 147.6565 & 17.2912, 16.0391, 16.7982 \\
$n=2000$ & 12.1743 & 193.6791 & 15.3551, 14.9715, 14.5999 \\
$n=3000$ & 11.4561 & 251.2571 & 14.3947, 13.8015, 14.0738 \\
$n=4000$ & 11.5022 & 310.8092 & 13.9249, 13.6388, 13.5541 \\
$n=5000$ & 10.8728 & 338.3008 & 13.2793, 12.9796, 13.3100 \\
  \hline
\end{tabular}
\end{center}
(The last column is based on three independent samples of $z \sim N(0,I_n)$.) Evidently,
the condition numbers do not increase with the dimension $n$ in the first and third columns, while
a near linear growth is observed in the second column.

We now illustrate (iii); moreover, we present experimental results suggesting
 that exact line search GD performs favorably compared to constant step size GD.
For the latter, we first show how to efficiently compute the line search function $p(t) := f(x+ td)$
by combining \eqref{eq:PhaseRetrievalQuartic} and \eqref{eq:PhaseRetrievalQuarticGrad}.
Write $A = [a_1, \cdots, a_m] \in \bR^{n \times m}$. We can compute the gradient descent direction,
the line search polynomial
$p(t)$, and the final `optimum' step size $s^\ast$
by computing a sequence of vectors and scalars in Algorithm~\ref{alg:ExactLineSearchPhaseRetrieval} below,
with the computational complexity of each step listed in curly braces.
\begin{algorithm}
    \caption{Computation of the exact line search step size $s^\ast$ for the (real) phase retrieval problem}
    \label{alg:ExactLineSearchPhaseRetrieval}
    \begin{algorithmic}[1]
    \STATE $z \leftarrow A^T x$ \;\; \COMMENT{\textcolor{blue}{$\approx 2mn$}}
    \STATE $\alpha \leftarrow  -y + z^2$ \;\; \COMMENT{$O(m)$}
    \STATE $d \leftarrow  -\nabla f(x) = -\frac{1}{m} A (\alpha \cdot z)$ \;\; \COMMENT{\textcolor{blue}{$\approx 2mn$}} 
    \STATE $w \leftarrow A^T d$ \;\; \COMMENT{\textcolor{blue}{$\approx 2mn$}}
    \STATE $\beta = 2 z \cdot w$ \;\; \COMMENT{$O(m)$}
    \STATE $\gamma = w^2$ \;\; \COMMENT{$O(m)$}
    \STATE $c_4\leftarrow \gamma^T \gamma$, $c_3\leftarrow 2 \beta^T \gamma$, $c_2 \leftarrow \beta^T\beta +2\alpha^T \gamma$,
           $c_1 \leftarrow 2\alpha^T \beta$, $c_0\leftarrow \alpha^T \alpha$ \;\; \COMMENT{$O(m)$}
    \STATE $s^\ast \leftarrow \argmin_{t\geq 0} p(t)$, \; $p(t)=f(x+td) = \sum_{i=0}^4 c_i t^i$ \;\; \COMMENT{$O(1)$}
    \end{algorithmic}
\end{algorithm}

In the algorithm description,  $u\cdot v$, for two vectors $u$, $v$ of the same length, stands for componentwise multiplication, and $v^2 :=v\cdot v$.
Notice that the exact step size $s^\ast$ in step 8 can be found in constant time by a root finder. (Since the degree is $4$, there is even a closed-formed formula for the roots.)
As we can see, the dominant steps are steps 1, 3 and 4. As only steps 1-3 are necessary for constant step size GD, we conclude that,
for the phase retrieval problem, \textbf{exact line search GD is about 50\% more expensive per iteration than constant step size GD}.

\newcommand{\parallelsum}{\mathbin{\!/\mkern-5mu/\!}}

\begin{figure}[ht]
\begin{center}
\begin{tabular}{c}
\includegraphics[height=5cm]{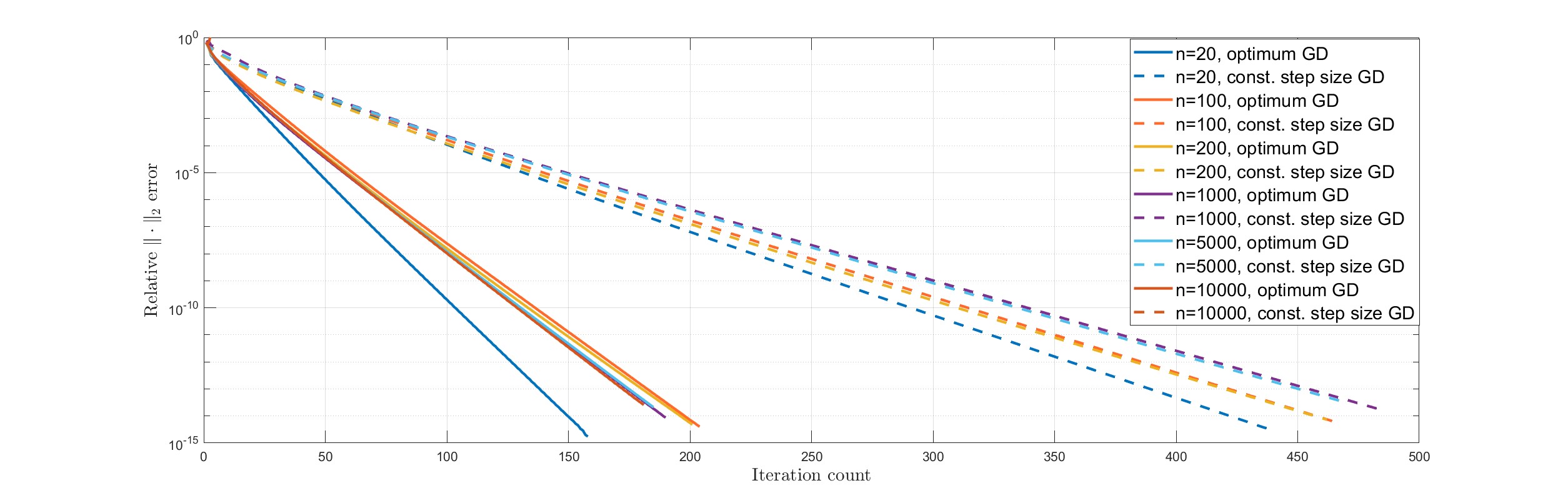}
\end{tabular}
\end{center}
\vspace{-.6cm}
\caption{ROC of constant step size GD vs optimum GD for the phase retrieval problem}
\label{fig:PhaseRetrievalROC}
\end{figure}

Figure~\ref{fig:PhaseRetrievalROC} shows the rates of convergence for gradient descent with constant step size $s=0.1$ (suggested in
\cite[Section 1.4]{MR4099988})
and exact line search
for $n=10, 100, 200, 1000, 5000, 10000$, $m=10n$, with the initial guess chosen by a so-called spectral initialization.\footnote{Under
the Gaussian design $a_1,\ldots,a_m \sim_{\rm i.i.d.} N(0,I_n)$, it is interesting to see that
$E[(a_i^T x)^2 a_i a_i^T] = I+2xx^T$. As the Rayleigh quotient $u^T (I+2xx^T) u /u^T u = 1+2 (u^Tx)^2/u^Tu$ is maximized when $u$ is parallel to $x$,
an educated initial guess for $x$ would be a suitably scaled dominant eigenvector of $\frac{1}{m} \sum_{r=1}^m y_r a_r a_r^T$, which can be computed
from the data $(y_r)_{r=1}^m$ and sensors $(a_r)_{r=1}^m$ using the power method. See, for example, \cite{MR4099988,MR3332993} for more details and references therein.}
As the plots show, for each signal size $n$
the ROC for exact line search GD is more than twice as fast as that of constant step size GD.

Not only is the speedup in ROC by exact line search outweighs
the 50\% increase in per-iteration cost,  the determination of step size is automatic and
requires no tuning -- arguably another practical advantage.
For each $n$, the ROC for exact line search GD in Figure~\ref{fig:PhaseRetrievalROC} is slightly faster than $O([(1-a)/(1+a)]^k)$ for
$a = {\rm cond}(\nabla^2 f(x^\ast))^{-1}$ -- the ROC attained by optimum GD as if the degree 4 objective has a constant Hessian $\nabla^2 f(x^\ast)$ (Theorem~\ref{thm:AkaikeBound}). 
Unsurprisingly, our experiments also suggest that
exact line search GD is more robust than its constant step size counterpart
for different choices of initial vectors.

Similar advantages for exact line search GD was observed for the matrix completion problem.
In ongoing work, we also found that Algorithm~\ref{alg:ExactLineSearchPhaseRetrieval} can be extended to the more realistic setting
of phase retrieval imaging from coded diffraction patterns \cite{CANDES2015277}, in which $A$ and $A^T$ become complex-valued and are implemented implicitly based on fast Fourier transforms. We shall report on these elsewhere.

\section{Invariance Properties of $\rho$ and $\GD$ and Akaike's $T$} \label{sec:Akaike}
Let $\bR^n_\ast$ be $\bR^n$ with the origin removed. For $x\in \bR^n$, define $|x|\in \bR^n$ by $|x|_i = |x_i|$.
Notice from \eqref{eq:rho} that, for a fixed $\lambda$, $\rho(\cdot,\lambda)$ is invariant under both scaling and sign-changes of the components, i.e.
\begin{align} \label{eq:invariance1}
\begin{split}
\rho(\alpha \mathcal{E} x, \lambda) = \rho(x,\lambda), \quad \forall x \in \bR^n_{\ast}, \; \alpha \neq 0, \;
\mathcal{E}={\rm diag}(\varepsilon_1,\ldots, \varepsilon_n), \; \varepsilon_i \in \{1,-1\}.
\end{split}
\end{align}
In other words, $\rho(x,\lambda)$ depends only on the equivalence class $[x]_\sim$, where
$x \sim y$ if $|x|/\|x\| = |y|/\|y\|$.
So, we overload/abuse notation and write
$$
\rho([x]_\sim,\lambda) := \rho(x,\lambda).
$$

By inspecting \eqref{eq:GD} one sees that
\begin{align} \label{eq:invariance2}
\begin{split}
\GD(\alpha \mathcal{E} x,\lambda) = \alpha \mathcal{E} \!\cdot\! \GD(x, \lambda),
\quad \forall x \in \bR^n_{\ast}, \; \alpha \neq 0, \;
\mathcal{E}={\rm diag}(\varepsilon_1,\ldots, \varepsilon_n), \; \varepsilon_i \in \{1,-1\}.
\end{split}
\end{align}
This means $[\GD(x)]_\sim$, when well-defined, depends only on $[x]_\sim$.
In other words, $\GD$ descends to a map
$[\GD]: {\rm dom}([\GD]) \subset \bR^n_{\ast}/\!\!\sim \rightarrow \bR^n_{\ast}/\!\!\sim.$
It can be shown that the natural domain of
$[\GD]$ is
\bea \label{eq:NaturalDomain}
\mathcal{M} :=  \big\{ [x]_\sim \in \bR^n_{\ast}/\!\!\sim \, : \mbox{$x$ is not an eigenvector of $A$} \big\};
\eea
also $[\GD](\mathcal{M}) \subset \mathcal{M}$.
It follows from the elegant variance-increasing property of Akaike's map $T$;
see \eqref{eq:AkaikeLemma2} below. But for clarity, we provide also a direct proof of this fact in Appendix~\ref{prop:Nonzero}.
So
\bea \label{eq:Tabstract}
[\GD]: \mathcal{M} \rightarrow \mathcal{M}
\eea
is a well-defined self-map.

Except when $n=2$ with $\lambda_1 > \lambda_2$,
$[\GD]$ does not extend continuously to the whole $\bR^n_{\ast}/\!\!\sim$; see below.

Thanks to these invariance properties, the ROC \eqref{eq:ROC_def} we aim to analyze can be rewritten as
\bea \label{eq:ROC_invar}
\rho^\ast(x^{(0)},\lambda) := \limsup_{k\goto \infty} \Bigg[ \prod_{j=0}^{k-1} \rho \Big( [\GD]^j \big( [x^{(0)}]_\sim \big), \lambda \Big) \Bigg]^{1/k}.
\eea
From this, we have:
\begin{proposition} \label{prop:ROC_T}
If $x^{(0)} \in \bR^n$ lies on one of the coordinate axis, then trivially $\rho^\ast(x^{(0)}, \lambda)=0$. Otherwise, $x^{(0)} \in \mathcal{M}$. In this case,
if 
 $X^{(\infty)} \subset \mathcal{M}$ is the set of cluster points of the orbit $\{[\GD]^j (x^{(0)}): j\geq 0\}$,
  then $$\rho^\ast(x^{(0)},\lambda) = \sup_{x \in X^{(\infty)}} \rho(x, \lambda).$$
\end{proposition}
As we shall see in Theorem~\ref{thm:AkaikeTheorem2}, for most $x^{(0)}$ the set of cluster points consists of two points in $\mathcal{M}$.

In other words, the ROC problem reduces to the limiting behavior of the discrete dynamical system defined by $\eqref{eq:Tabstract}$.

\gap
\noindent{\bf Akaike's map $T$.}
While Akaike, a statistician, did not use jargons such as `invariance', `parametrization' or `dynamical systems' in his paper, the map
$T$ introduced in \cite{Akaike-1959} is the representation of $[\GD]$
under the identification of $[x]_\sim$ with
\bea \label{eq:sigma}
\sigma([x]_\sim) := \Big[ \lambda_1^2 x_1^2, \ldots, \lambda_n^2 x_n^2 \Big]^T / \sum_{j} \lambda_j^2 x_j^2 \in \triangle_n := \Big\{ p \in \bR^n: \sum_j p_j = 1, p_j \geq 0 \Big\}.
\eea
In above, $\triangle_n$, or simply $\triangle$, is usually called the standard simplex, or the probability simplex as Akaike would prefer.
One can verify that $\sigma: \bR^n_{\ast}/\!\!\sim \goto \triangle$ is a well-defined bijection and hence
\bea \label{eq:sigmaInv}
\sigma^{-1}: \triangle \goto \bR^n_{\ast}/\!\!\sim, \quad p \mapsto [x]_\sim, \; x_j = \frac{\sqrt{p_j}}{\lambda_j}
\eea
 may be viewed as a
parametrization of the quotient space $\bR^n_{\ast}/\!\!\sim$.
(Strictly speaking, the map $\sigma^{-1}$ is not a parametrization. As a manifold,
$\bR^n_{\ast}/\!\!\sim$ is $(n-1)$-dimensional, which means it deserves a parametrization with $n-1$ parameters.
But, of course, we can identify any $p \in \triangle$ with $[s_1,\ldots,s_{n-1}]^T$ by
$p = [s_1, \ldots, s_{n-1}, 1-\sum_{i=1}^{n-1} s_i]^T$; we shall do exactly this in \eqref{eq:Remove1dof} while proving Theorem~\ref{thm:AkaikeBound}.)

We now derive a formula for $T:= \sigma \circ [\GD] \circ \sigma^{-1}$:
By \eqref{eq:GD} and \eqref{eq:sigmaInv}, $[\GD](\sigma^{-1}(p))$ has a representor $y\in \bR^n_{\ast}$ with
\begin{align*}
y_i &= \frac{\sqrt{p_j}}{\lambda_j} - \frac{\sum_j \lambda_j^2 p_j/\lambda_j^2}{\sum_j \lambda_j^3 p_j/\lambda_j^2} \lambda_i \frac{\sqrt{p_j}}{\lambda_j}
   = \sqrt{p_i} \Big[ \frac{1}{\lambda_i} - \frac{1}{\sum_j \lambda_j p_j} \Big]
 = \sqrt{p_i} \, \frac{\overline{\lambda}(p)-\lambda_i}{\lambda_i \overline{\lambda}(p)},
\end{align*}
where $\overline{\lambda}(p) := \sum_j \lambda_j p_j$.
Consequently,
\bea \label{eq:T}
T(p)_i = \lambda_i^2 \Big( \sqrt{p_i} \, \frac{\overline{\lambda}(p)-\lambda_i}{\lambda_i \overline{\lambda}(p)} \Big)^2 / \sum_j \lambda_j^2 \Big( \sqrt{p_j} \, \frac{\overline{\lambda}(p)-\lambda_j}{\lambda_j \overline{\lambda}(p)} \Big)^2
= \frac{p_i (\overline{\lambda}(p)-\lambda_i)^2}{\sum_j p_j (\overline{\lambda}(p)-\lambda_j)^2}.
\eea
The last expression is Akaike's map $T$
defined in \cite[\S1]{Akaike-1959}.
Under the distinct eigenvalues assumption (see below), $T(p)$ is well-defined for any $p$ in
\bea \label{eq:naturaldomain}
{\rm dom}(T) = \triangle_n \backslash \{e_1,\ldots,e_n\},
\eea
i.e. the standard simplex with the standard basis of $\bR^n$ removed.
Also, $T$ is continuous on its domain.
By \eqref{eq:T2d} below, when $n=2$, $T$ extends continuously to $\triangle_2$.
But for $n\geq 3$ it does not extend
continuously to any $e_i$; for example, if $n=3$ and $i=2$, then (assuming $\lambda_1>\lambda_2>\lambda_3$),
\bea \label{eq:diagonal_n3}
T([\epsilon,1-\epsilon,0]^T)=[1-\epsilon,\epsilon,0]^T \quad \mbox{and} \quad
T([0,1-\epsilon,\epsilon]^T)=[0,\epsilon,1-\epsilon]^T.
\eea
This follows from the $n=2$ case of Proposition~\ref{prop:Tinvariance} and the following diagonal property of $T$.

\gap
\noindent
{\bf Diagonal property.} Thanks to the matrix $A$ being diagonal, $T$ is invariant under
$\triangle_J := \{ p \in \triangle_n : p_i = 0, \; \forall i \notin J \}$ for any
$J \subset \{1,\ldots,n\}$. Notice the correspondence between $\triangle_J$ and $\triangle_{|J|}$
via the projection $\lambda \mapsto \lambda_J := (\lambda_i)_{i \in J}$.
If we write $T_\lambda$ to signify the dependence of $T$ on $\lambda$,
then under this correspondence $T|_{\triangle_J}$ is simply $T_{\lambda_J}$ acting on $\triangle_{|J|}$.

This obvious property will be useful for our proof later; for now, see \eqref{eq:diagonal_n3}
above for an application of the property.

\gap
\noindent
{\bf Distinct eigenvalues assumption.} It causes no loss of generality to assume that the eigenvalues $\lambda_i$ are distinct:
if $\hat{\lambda} = [\hat{\lambda}_i, \ldots, \hat{\lambda}_m]^T$ consists of the distinct eigenvalues of $A$, then $A = {\rm diag}(\hat{\lambda}_i I, \ldots, \hat{\lambda}_m I)$,
where each $I$ stands for an identity matrix of the appropriate size. Accordingly, each
initial vector $x^{(0)}$ can be written in block form $[\mathbf{x}^{(0)}_1, \ldots, \mathbf{x}^{(0)}_m]^T$.
It is easy to check that if we apply the optimum GD method to $\hat{f}(\hat{x}) = \frac{1}{2} \hat{x}^T \hat{A} \hat{x}$
with $\hat{A}={\rm diag}(\hat{\lambda}_i, \ldots, \hat{\lambda}_m)$ and
 initial vector $\hat{x}^{(0)} = \big[\|\mathbf{x}^{(0)}_1\|_2, \ldots, \|\mathbf{x}^{(0)}_m\|_2\big]^T$, then the ROC of the reduced system is identical to that of the original.
Moreover, the map 
$$P: \bS^{n-1}_+ \goto  \bS^{m-1}_+, \quad [\mathbf{x}_1, \ldots, \mathbf{x}_m]^T \mapsto
\big[\|\mathbf{x}_1\|_2, \ldots, \|\mathbf{x}_m\|_2\big]^T,$$
is a submersion and hence
has the property that $P^{-1}(N)$ is a null set in $\bS^{n-1}$ for any null set $N$ in $\bS^{m-1}_+$; see Lemma~\ref{lemma:Null1}.
Therefore, it suffices to prove Theorem~\ref{thm:Main} and Theorem~\ref{thm:AkaikeBound} under the distinct eigenvalues assumption.

So from now on we make the blanket assumption that $\lambda_1>\cdots>\lambda_n>0$.

\gap
\noindent
{\bf Probabilistic interpretation.}
Akaike's $\lambda$-dependent parametrization \eqref{eq:sigma}-\eqref{eq:sigmaInv} does not only give $[\GD]$ the simple representation $T$
\eqref{eq:T}, the map $T$ also has an interesting probabilistic interpretation:
if we think of $p$ as a probability distribution associated to the values in $\lambda$, then $\overline{\lambda}(p)$ is the mean of the resulted random variable,
the expression in the dominator of the definition of $T$, i.e. $\sum_j p_j (\overline{\lambda}(p)-\lambda_j)^2$, is the variance of the random variable.
What, then, does the map $T$ do to $p$? It produces a new probability distribution, namely $T(p)$, for $\lambda$.
The definition of $T$ in \eqref{eq:T} suggests that $T(p)_i$ will be bigger if $\lambda_i$ is far from the mean $\overline{\lambda}(p)$, so the map
polarizes the 
probabilities to the extremal values $\lambda_1$ and $\lambda_n$.
This also suggests that the map $T$ tends to increase variance.
Akaike proved that it is indeed the case in \cite[Lemma 2]{Akaike-1959}:
Using his notation $\overline{f(\lambda)}(p) := \sum_{i=1}^n f(\lambda_i) p_i$, the result can be expressed as
\bea \label{eq:AkaikeLemma2}
\overline{(\lambda - \overline{\lambda}(T(p)) )^2}(T(p)) \geq \overline{(\lambda - \overline{\lambda}(p) )^2}(p).
\eea
This variance increasing property is a key to Akaike's proof of \eqref{eq:AkaikeTheorem12} below. As an
immediate application, notice that by \eqref{eq:naturaldomain}
$p \in {\rm dom}(T)$ is equivalent to saying that the random variable with probability $p_i$ attached to the value $\lambda_i$ has a positive variance.
Therefore \eqref{eq:AkaikeLemma2} implies that if $p \in {\rm dom}(T)$, then $T(p) \in {\rm dom}(T)$ also, and so $T^k(p)$ is well-defined for all $k \geq 0$.
This shows that
$$
T: \triangle_n \backslash \{e_1,\ldots,e_n\} \goto \triangle_n \backslash \{e_1,\ldots,e_n\}
$$
is a well-defined self-map. This is, again, the abstractly defined map \eqref{eq:Tabstract} parameterized in a very special way.

\gap
The following fact is instrumental for our proof of the main result; it not explicitly stated in \cite{Akaike-1959}.
\begin{proposition}[Independence of $\lambda_1$ and $\lambda_n$] \label{prop:Tinvariance}
The map $T$ depends only on $\alpha_i\in (0,1)$, $i=2,\ldots, n-1$, defined by
\bea \label{eq:alpha}
\lambda_i = \alpha_i \lambda_1 + (1-\alpha_i) \lambda_n.
\eea
In particular, when $n=2$, $T$ is independent of $\lambda$; in fact,
\bea \label{eq:T2d}
T([s,1-s]^T) = [1-s,s]^T.
\eea
\end{proposition}
\pf
Write $\alpha_1:=1$, $\alpha_n:=0$.
The map $T$ can be expressed as
\bea \label{eq:T_in_alphas}
T(p)_i = \frac{p_i \big[\sum_\ell p_\ell (\lambda_\ell-\lambda_i)\big]^2}{\sum_j p_j \big[\sum_\ell p_\ell (\lambda_\ell-\lambda_j)\big]^2}
= \frac{p_i \big[\sum_\ell p_\ell (\alpha_\ell-\alpha_i)\big]^2}{\sum_j p_j \big[\sum_\ell p_\ell (\alpha_\ell-\alpha_j)\big]^2}.
\eea
The first equality is due to $\sum_i p_i=1$, whereas the second equality is due to
$
\lambda_\ell-\lambda_i = \alpha_\ell \lambda_1 + (1-\alpha_\ell) \lambda_n - \alpha_i \lambda_1 - (1-\alpha_i) \lambda_n
=(\alpha_\ell - \alpha_i) (\lambda_1-\lambda_n).
$
\eop

While the proof is trivial, it is not clear if Akaike was aware of the first part of the above proposition. It tells us that
the condition number of $A$ does not play a role in the dynamics of $T$. However, he must be aware of the second part
of the proposition, as he proved the following
 nontrivial generalization of \eqref{eq:T2d} in higher dimensions.
When the dimension is higher than 2, $T$ is no longer an involution, but yet $T$ resembles \eqref{eq:T2d}
in the following way:

\begin{theorem}{\cite[Theorem 2]{Akaike-1959}} \label{thm:AkaikeTheorem2}
For any $p^{(0)} \in {\rm interior}(\Delta_n)$ (i.e. $p^{(0)}_i>0$, $\forall i$), there exists $s \in (0,1)$ such that
\bea \label{eq:AkaikeTheorem12}
p^{(\infty)} := \lim_{k \goto \infty} T^{2k}(p^{(0)}) = [1-s,0,\ldots,0,s]^T \; \mbox{ and } \;
p^{\ast(\infty)} := \lim_{k \goto \infty} T^{2k+1}(p^{(0)}) = [s,0,\ldots,0,1-s]^T.
\eea
\end{theorem}

This theorem 
says that the dynamical system defined by $T$ polarizes the probabilities to
the two extremal eigenvalues. Combining with Proposition~\ref{prop:ROC_T}, the ROC problem reduces to understanding the limit probability $s$ in \eqref{eq:AkaikeTheorem12}.

This makes part of the general problem two-dimensional; so
we begin in Section~\ref{sec:2D} by analyzing the $n=2$ case and proving Theorem~\ref{thm:Main}(i).
The proof of Theorem~\ref{thm:AkaikeBound}, which addresses the $n>2$ case, to be presented in Section~\ref{sec:TheoremAkaikeBound}, utilizes some of the intermediate calculations in the 2-D case, but the key technical contribution is a refinement of Akaike's Theorem~\ref{thm:AkaikeTheorem2}, which we briefly illustrate in the next subsection.

\subsection{An overview of the proof of Theorem~\ref{thm:AkaikeBound}}
\begin{figure}[ht]
\begin{center}
\begin{tabular}{c}
\includegraphics[height=5cm]{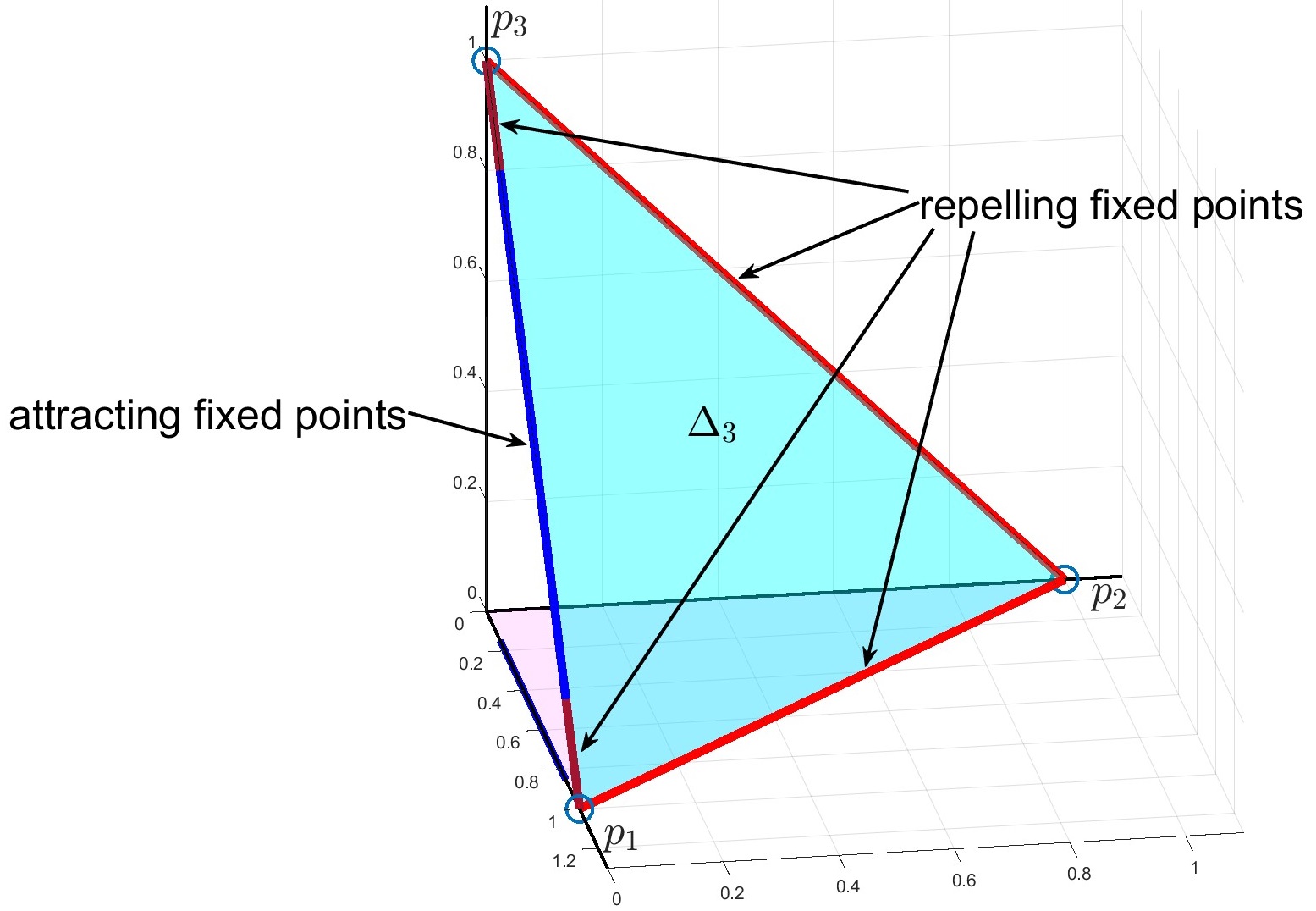}
\end{tabular}
\end{center}
\vspace{-.6cm}
\caption{In dimension $n=3$, every point on the boundary of $\Delta_3\backslash \{e_1,e_2,e_3\}$ is a fixed point of $T^2=T\circ T$.
Every point in the interior of $\Delta_3$ is attracted by $T^2$ to some point in the blue line segment.
The analysis is done by writing $T$ as a map on the projected simplex of $\Delta$ on the $p_1$-$p_2$ plane, colored in light purple; see \eqref{eq:Remove1dof}.}
\label{fig:AkaikeMap_3}
\end{figure}

When $n=2$ the map $T$ is simply given by \eqref{eq:T2d} and its dynamics is trivial.
So as far as the dynamics of $T$ is concerned, the challenge begins with
$n=3$, which is also the lowest dimension that allows for an intermediate eigenvalue.
In this case, $T$ maps the simplex $\Delta_3$ with the corners removed, as shown in Figure~\ref{fig:AkaikeMap_3}, to itself.
By \eqref{eq:T_in_alphas}, the map depends only on $\alpha_2 \in (0,1)$, defined by $\lambda_2 = \alpha_2 \lambda_1 + (1-\alpha_2) \lambda_3$).

When restricted to each of the three edges of the triangle $\Delta_3$, $T$ behaves exactly like its 2-D counterpart \eqref{eq:T2d}
(recall \eqref{eq:diagonal_n3}), so
every point on the boundary of ${\rm dom}(T)$ is a fixed point of $T^2 = T\circ T$.
According to Theorem~\ref{thm:AkaikeTheorem2}, every point in the interior of ${\rm dom}(T)$ is attracted by $T^2$ to
a point of the form $[s,0,1-s]^T$, $s \in (0,1)$, i.e. $\lim_{k\goto\infty} T^{2k}(p) =[s,0,1-s]^T$.

The key mathematical finding leading to a proof of Theorem~\ref{thm:AkaikeTheorem2}
-- and a key difference between the $n=2$ and $n=3$ cases --
is that
the limit probability $s$ in \eqref{eq:AkaikeTheorem12} above lies on an interval strictly smaller than $(0,1)$.
When $\alpha_2 = 1/2$, the interval is
\bea \label{eq:I}
I = [1/2-1/(2\sqrt{2}), 1/2+1/(2\sqrt{2})].
\eea
The interval
is illustrated by the blue line segment in Figure~\ref{fig:AkaikeMap_3}.

More generally:
\begin{itemize}
\item The interval $I$ depends on the $\alpha_i$ closest to $1/2$
according to \eqref{eq:AttractingInterval} in Lemma~\ref{lemma:Spectrum}.
\item When $n\geq 4$, the set of points $p^{(0)}$ in ${\rm interior}(\Delta_n)$ that satisfy \eqref{eq:AkaikeTheorem12} with $s \notin I$ is \textbf{non-empty but null},
i.e. \eqref{eq:AkaikeTheorem12} holds with $s \in I$ for \emph{almost all $p^{(0)}\in {\rm interior}(\Delta_n)$}.
\end{itemize}
This is Claim IV in the proof of Theorem~\ref{thm:AkaikeBound}. While its rigorous proof is technical and requires the theorem of center and stable manifolds (Theorem~\ref{thm:StableManifolds}), along with other arguments, it is relatively easy to see how the interval $I$ comes about:
it follows from calculating the spectrum of the Jacobian of $T^2$ at the fixed points $[s,0,\ldots,0,1-s]^T$ and
see what values of $s$ give ``attracting fixed points", i.e. the spectrum consists of eigenvalues with
modulus no bigger than 1. The interval $I$ is precisely the set of such values of $s$.

To the surprise of the author, the interval $I$, derived using the aforementioned dynamical system approach, corresponds exactly to
Akaike's original ROC bound (right-hand side of \eqref{eq:LowerBound}), which was derived using a tricky but seemingly
haphazard argument in \cite[\S2]{Akaike-1959}, resulting also in a difficult-to-interpret condition on $x^{(0)}$ for the lower bound in \eqref{eq:LowerBound} to hold true.
His argument makes it hard to see whether the lower bound is tight, while his condition on $x^{(0)}$ makes it hard to
see if it is satisfied by almost all $x^{(0)}$ in $\Delta_n$.

\section{Analysis in 2-D and Proof of Theorem~\ref{thm:Main}(i)}  \label{sec:2D}
When $n=2$, we may represent a vector in $\triangle$ by $[1-s,s]^T$ with $s\in[0,1]$.
Recall \eqref{eq:T2d} and note that
$\rho$ depends on $\lambda$ only through $a:= \lambda_2/\lambda_1 = {\rm cond}(A)^{-1}$.
So, we may represent $T$ and $\rho$ in the parameter $s$ and the quantity $a$ as:
\bea \label{eq:t_rho}
t(s)=1-s, \quad \overline{\rho}^2(s,a)=1-(1-s+sa)^{-1} (1-s+s a^{-1})^{-1}.
\eea
So  $\overline{\rho}({t}(s),a) = \overline{\rho}(s,a)$, and
the otherwise difficult-to-analyze ROC $\rho^\ast$ (\eqref{eq:ROC_def}) is determined simply by
\bea \label{eq:rhostar_rho}
\rho^\ast(x^{(0)},\lambda) = \rho(x^{(0)},\lambda).
\eea

By \eqref{eq:t_rho}, the value $s\in [0,1]$ that maximizes $\overline{\rho}$ 
is given by $s_{\rm max}=1/2$, with
maximum value $(1-a)/(1+a)$ ($=(\lambda_1-\lambda_2)/(\lambda_1+\lambda_2)$).
This means, by \eqref{eq:sigma}-\eqref{eq:sigmaInv}, $\rho^\ast(x^{(0)},\lambda) = \rho(x^{(0)},\lambda) = (1-a)/(1+a)$ for any $x^{(0)}\in \bR_\ast^2$ with
$$
|x^{(0)}_1| = \lambda_2/\lambda_1 |x^{(0)}_2|.
$$
This observation also shows that the worst case bound \eqref{eq:Worst} is tight in any dimension $n$:
\begin{proposition} \label{prop:SharpestOneStep}
There exists initial vector $x^{(0)} \in \bR^n$ such that
equality holds in \eqref{eq:Worst} for any iteration $k$.
\end{proposition}
\pf
We have proved the claim for $n=2$.
For a general dimension $n\geq 2$, observe that if the initial vector lies on the $x_1$-$x_n$ plane, then -- thanks to diagonalization --
GD behaves exactly the same as in 2-D, with $A={\rm diag}(\lambda_1,\lambda_2,\ldots,\lambda_n)$
 replaced by $A = {\rm diag}(\lambda_1,\lambda_n)$.
That is, if we
choose $x^{(0)} \in \bR^n$ so that $|x^{(0)}_1| = \lambda_n/\lambda_1 |x^{(0)}_n|$
and $x^{(0)}_i = 0$ for $2 \leq i \leq n-1$, then equality holds in \eqref{eq:Worst} for any $k$.
\eop

\begin{figure}[ht]
\begin{center}
\begin{tabular}{cc}
\includegraphics[height=6cm]{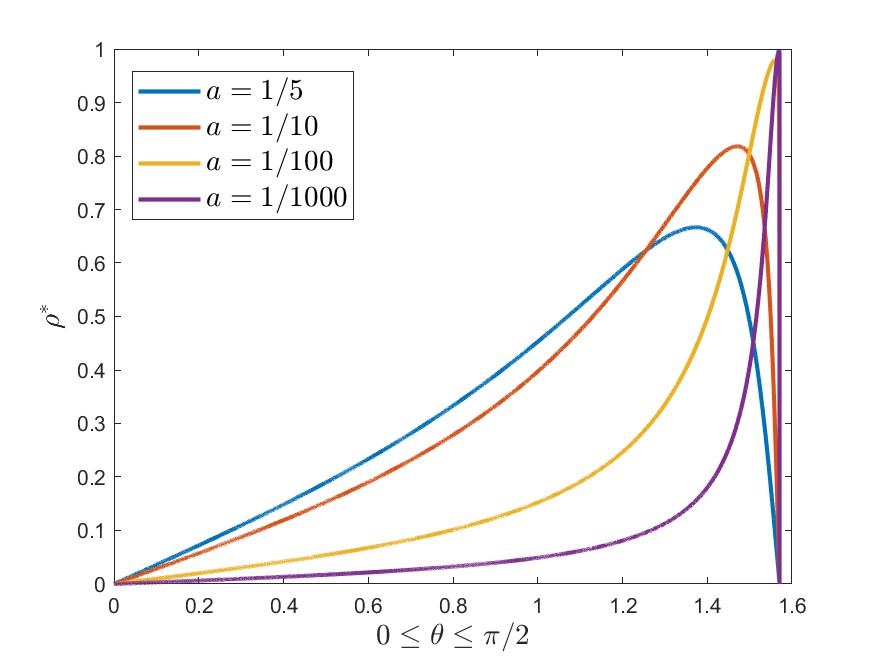}
&
\includegraphics[height=6cm]{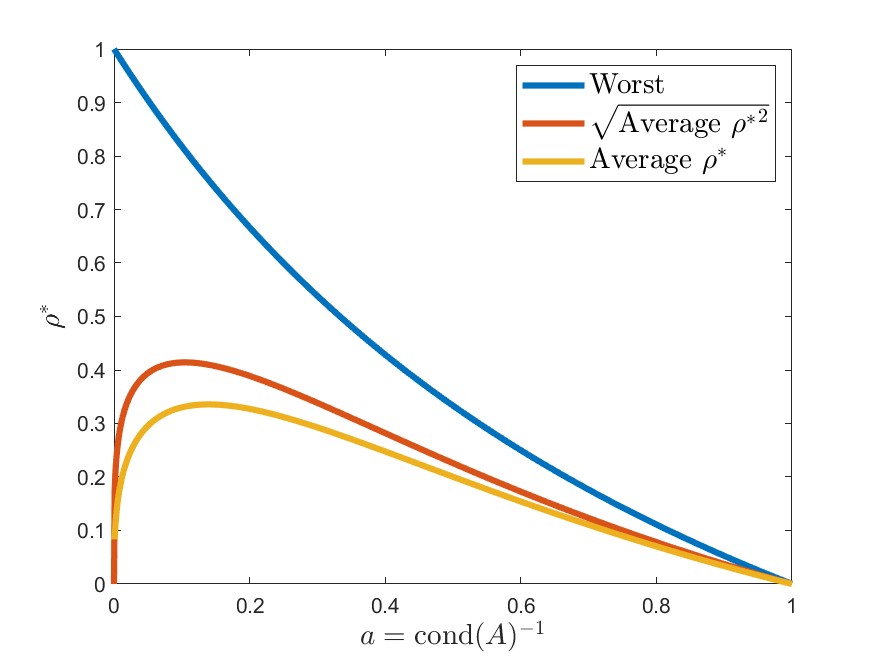}
\end{tabular}
\end{center}
\vspace{-.6cm}
\caption{Left: Plots of $\theta$ versus $\rho^\ast([\cos(\theta),\sin(\theta)]^T,[1,a]^T)$ for various values of $a$.
Observe the convergence in \eqref{eq:NonUniform} being non-uniform in $\theta$.
Right: The worst, average and the square root of the average square rate of convergence as a function of $a$.
The average rate of convergence is computed using numerical integration, while the other two curves are given by the
closed-form expressions \eqref{eq:Worst2D} and \eqref{eq:closed}.}
\label{fig:WorstvsAverage}
\end{figure}

\gap
\noindent
{\bf Proof of Theorem~\ref{thm:Main}(i).}
For any non-zero vector $x$ in $\bR^2$, $[x]_\sim$ can be identified by $[\cos(\theta),\sin(\theta)]^T$ with $\theta \in [0,\pi/2]$, which,
by \eqref{eq:sigma}-\eqref{eq:sigmaInv}, is identified with $[1-s,s]^T \in \Delta_2$ where
\bea \label{eq:theta_s}
s=\frac{a^2\sin^2(\theta)}{\cos^2(\theta)+a^2\sin^2(\theta)}.
\eea
The average ROC (recall Definition~\ref{eq:AveROC}) is given by
\bea \label{eq:Ave2D}
\mbox{Average ROC} = (\pi/2)^{-1} \int_0^{\pi/2} \rho([\cos(\theta), \sin(\theta)]^T,[1,a]) \,d\theta.
\eea
Note that
\bea \label{eq:rho2}
\rho([\cos(\theta), \sin(\theta)]^T,[1,a]^T) = \overline{\rho}\Big(\frac{a^2\sin^2(\theta)}{\cos^2(\theta)+a^2\sin^2(\theta)},a\Big).
\eea
On the one hand, we have
\bea \label{eq:Worst2D}
\max_{\theta \in [0,\pi/2]} \rho([\cos(\theta), \sin(\theta)]^T,[1,a]^T) = \overline{\rho}(1/2,a) = (1-a)/(1+a),
\eea
so
$\lim_{a \goto 0^+} \max_{s\in[0,1]} \overline{\rho}(s,a) = 1$.
On the other hand, by \eqref{eq:rho2} and \eqref{eq:t_rho} one can verify that
\bea \label{eq:NonUniform}
 \lim_{a \goto 0^+} \rho([\cos(\theta), \sin(\theta)]^T,[1,a]^T) = 0, \;\; \forall\, \theta \in [0,\pi/2].
\eea
See Figure~\ref{fig:WorstvsAverage} (left panel) illustrating the non-uniform convergence.
Since $\rho([\cos(\theta), \sin(\theta)]^T,[1,a]^T) \leq 1$,
by the dominated convergence theorem,
$$
\lim_{a\goto 0^+} \int_0^{\pi/2} \rho([\cos(\theta), \sin(\theta)]^T,[1,a]^T) \,d\theta =
\int_0^{\pi/2} \lim_{a\goto 0^+} \rho([\cos(\theta), \sin(\theta)]^T,[1,a]^T) \,d\theta = 0.
$$
This proves the first part of Theorem~\ref{thm:Main}. \eop

\noindent
{\bf An alternate proof.}
While the average rate of convergence $(\pi/2)^{-1} \int_0^{\pi/2} \rho([\cos(t),\sin(t)]^T,a) \,d\theta$
does not seem to have a closed-form expression, the average \emph{square} rate of convergence can be
expressed in closed-form:
\bea \label{eq:closed}
(\pi/2)^{-1} \int_0^{\pi/2} \overline{\rho}^2\Big(\frac{a^2\sin^2(\theta)}{\cos^2(\theta)+a^2\sin^2(\theta)},a\Big) \,d\theta =
\frac{ \sqrt{a} (1-\sqrt{a})^2}{(1+a)(1-\sqrt{a}+a)}. 
\eea
By Jensen's inequality,
\bea \label{eq:Jensen}
\big[ (\pi/2)^{-1} \int_0^{\pi/2} \rho([\cos(\theta),\sin(\theta)]^T,\lambda) \,d\theta \big]^2 \leq (\pi/2)^{-1} \int_0^{\pi/2} \rho^2([\cos(\theta),\sin(\theta)]^T,\lambda) \,d\theta.
\eea
Since the right-hand side of \eqref{eq:closed} ($=$ the r.h.s. of \eqref{eq:Jensen}) goes to zero as $a$ approaches 0, so does the left-hand side of \eqref{eq:Jensen} and
hence also the average ROC. \qed

See Figure~\ref{fig:WorstvsAverage} (right panel, yellow curve) for a plot of how the average ROC varies with ${\rm cond}(A)$: as ${\rm cond}(A)$ increases from 1, the average
ROC does deteriorate -- as most textbooks would suggest -- but only up to a certain point, after that the average ROC does not only improve, but gets arbitrarily fast,
as $A$ gets more and more ill-conditioned -- quite the opposite of what most textbooks may suggest. See also Appendix~\ref{AppendixC}.

\section{Proof of Theorem~\ref{thm:AkaikeBound} and Theorem~\ref{thm:Main}(ii)} \label{sec:TheoremAkaikeBound}
Let $n\geq 3$.
Denote by $\Theta$ the map
\bea \label{eq:Remove1dof}
[p_i]_{1\leq i<n} \mapsto \big[ T_i(p_1,\ldots,p_{n-1},1-\mbox{$\sum_{j=1}^{n-1} p_j$}) \big]_{1\leq i<n}.
\eea
Its domain, denoted by ${\rm dom}(\Theta)$, is the
simplex $\Lambda := \big\{[p_1,\ldots,p_{n-1}]^T: p_j\geq 0, \; 0 \leq \sum p_j \leq 1 \big\}$ with its vertices removed.
Notwithstanding, $\Theta$ can be
be smoothly extended to some open set of $\bR^{n-1}$ containing ${\rm dom}(\Theta)$.

The 2-periodic points $[s,0,\ldots,0,1-s]^T$ of $T$ according to $\eqref{eq:AkaikeTheorem12}$ correspond to the fixed points
$[s,0,\ldots,0]^T$, which we denote more compactly by
$\big[\!\begin{smallmatrix} s \\ 0 \end{smallmatrix}\!\big]$, of $\Theta^2 = \Theta \circ \Theta$.

The map $\sigma$ defined by \eqref{eq:sigma}-\eqref{eq:sigmaInv} induces smooth one-to-one correspondences between $\bS^{n-1}_+$, $\triangle$ and $\Lambda$.
For any $x^{(0)} \in \bS^{n-1}_+$, let $p^{(0)} \in \triangle$ be the corresponding probability vector and denote by
$s(x^{(0)}) \in (0,1)$
the corresponding \textbf{limit probability} according to $\eqref{eq:AkaikeTheorem12}$.
Theorem~\ref{thm:AkaikeTheorem2}, together with \eqref{eq:t_rho}, imply that
\begin{align}
\nonumber
\rho^\ast(x^{(0)},\lambda) & = \sqrt{1-\big( 1-s+sa \big)^{-1} \big( 1-s+s a^{-1} \big)^{-1}}, \quad s = s(x^{(0)}) \\
& = \frac{1-a}{\sqrt{(1+a)^2 + a(c - c^{-1})^2}} \quad \mbox{ if we write $s(x^{(0)})=1/(1+c^2)$.}  \label{eq:ROC}
\end{align}

\begin{figure}[ht]
\begin{center}
\begin{tabular}{cc}
\includegraphics[height=6cm]{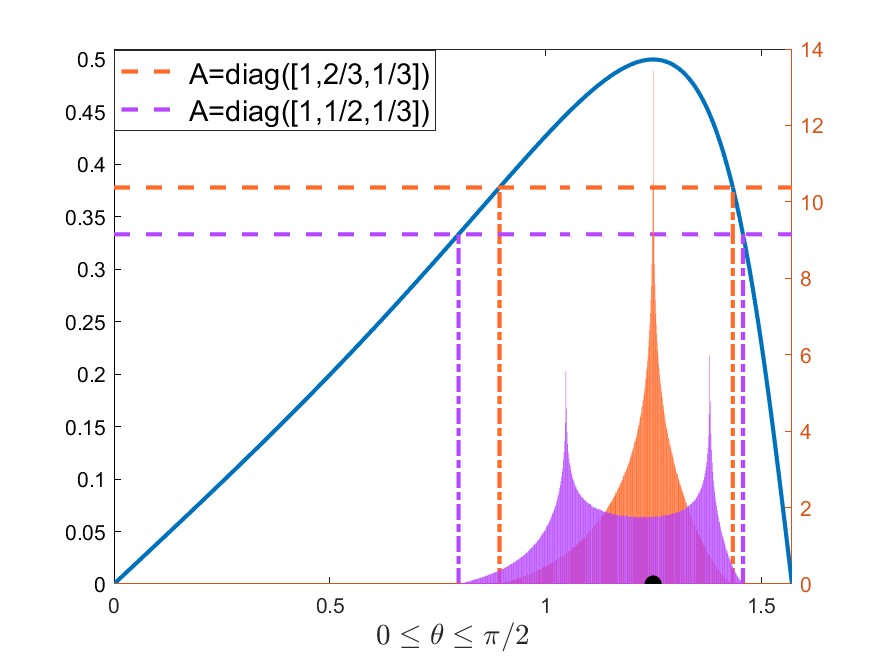}
&
\includegraphics[height=6cm]{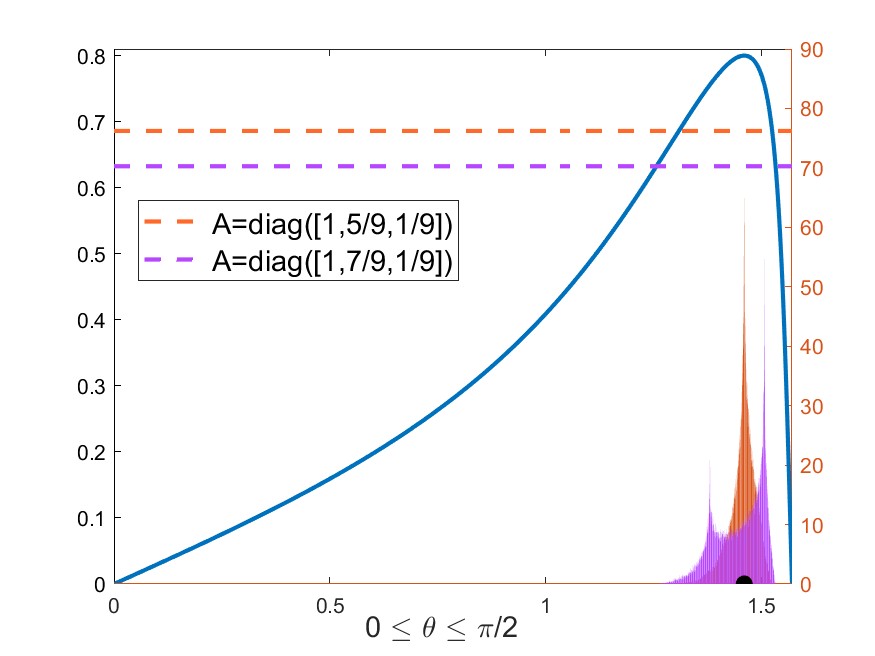}
\end{tabular}
\end{center}
\vspace{-.6cm}
\caption{Distribution of the limit angle $\theta$, estimated from $10^7$ initial vectors $x^{(0)}$ sampled from the uniform distribution on the unit 2-sphere in 3-D. The distribution is supported on a sub-interval $J$ of $[0,\pi/2]$, where $J$ is defined by \eqref{eq:Interval_J}.
The horizontal lines show the essential best ROC in \eqref{eq:LowerBound}. The left vertical axis is for ROC, while the right vertical axis is for probability density.
The black dot corresponds to the angle $\tan^{-1}(a^{-1})$, which yields the slowest ROC $(1-a)/(1+a)$.}
\label{fig:Figure4}
\end{figure}

A computation shows that for any $s \in (0,1)$, the Jacobian matrix of $\Theta$ at $\big[\!\begin{smallmatrix} s \\ 0 \end{smallmatrix}\!\big]$
is:
\bea \label{eq:DTheta_FixedPoint}
D\Theta|_{\big[\!\begin{smallmatrix} s \\ 0 \end{smallmatrix}\!\big]} :=
\frac{\partial(\Theta_1,\ldots,\Theta_{n-1})}{\partial(p_1,\ldots,p_{n-1})}\Big|_{\big[\!\begin{smallmatrix} s \\ 0 \end{smallmatrix}\!\big]} =
\begin{bmatrix}
  -1 & -\frac{\alpha_2^2}{s} & \cdots & \cdots & -\frac{\alpha_{n-1}^2}{s} \\
  0 & \frac{(\alpha_2 - s)^2}{s(1-s)} & 0 & \cdots & 0 \\
  \vdots & \ddots & \ddots & \ddots &  \vdots \\
  \vdots & \ddots & \ddots & \ddots &  0 \\
  0 & \cdots &  \cdots & 0 & \frac{(\alpha_{n-1} - s)^2}{s(1-s)} \\
\end{bmatrix},
\eea
where $\alpha_i$ is defined as in Proposition~\ref{prop:Tinvariance}.
Since $\Theta(\big[\!\begin{smallmatrix} s \\ 0 \end{smallmatrix}\!\big]) = \big[\!\begin{smallmatrix} 1-s \\ 0 \end{smallmatrix}\!\big]$,
the Jacobian matrix of $\Theta^2$ at its fixed point $\big[\!\begin{smallmatrix} s \\ 0 \end{smallmatrix}\!\big]$ is given by
$$
D \Theta^2|_{\big[\!\begin{smallmatrix} s \\ 0 \end{smallmatrix}\!\big]} =
D\Theta|_{\big[\!\begin{smallmatrix} 1-s \\ 0 \end{smallmatrix}\!\big]} \cdot
D\Theta|_{\big[\!\begin{smallmatrix} s \\ 0 \end{smallmatrix}\!\big]};
$$
its eigenvalues are 1 and
\bea \label{eq:SpectrumPhi}
\mu_i(s) := \frac{(\alpha_i - s)^2 (\alpha_i - (1 - s))^2}{s^2(1-s)^2} = \left(\frac{s(1-s)-\alpha_i(1-\alpha_i)}{s(1-s)} \right)^2, \;\; i=2,\ldots,n-1.
\eea
Each of the last $n-2$ eigenvalues, namely $\mu_i(s)$, is less than or equal to 1 if and only if $s (1-s) \geq \frac{1}{2} \alpha_i(1-\alpha_i)$.
Consequently,
we have the following:
\begin{lemma} \label{lemma:Spectrum} The spectrum of $D\Theta^2|_{\big[\!\begin{smallmatrix} s \\ 0 \end{smallmatrix}\!\big]}$ has at least one eigenvalue larger than one
iff $s \in (-1,1) \backslash I$, where
\bea \label{eq:AttractingInterval}
I := \Big\{ s: |s-1/2| \leq \frac{1}{2} \sqrt{1-2\alpha_{i^\ast}(1-\alpha_{i^\ast})} \Big\},
\;\;\;
i^\ast \in \argmin_{1<i<n} |\alpha_{i} - 1/2| \;
(=\argmin _{1<i<n} |\lambda_{i} - (\lambda_1+\lambda_n)/2|.)
\eea
\end{lemma}
Next, observe that:

\gap
\noindent
{\bf Claim I}:
If $x^{(0)} \in \bS^{n-1}_+$ satisfies the upper bound on $|s(x^{(0)})-1/2|$ in the definition of $I$, then
the ROC $\rho^\ast(x^{(0)},\lambda)$ satisfies the lower bound in \eqref{eq:LowerBound}.
This can be checked by a straightforward computation using the expression of $\rho^\ast(x^{(0)},\lambda)$ in \eqref{eq:ROC}.
(The correspondence of the interval $I$ and the lower bound \eqref{eq:LowerBound} is illustrated in Figure~\ref{fig:Figure4}.)

\gap
\noindent
{\bf Claim II}:  The correspondence between $\bS^{n-1}_+$ and $\Lambda$, induced by \eqref{eq:sigma}-\eqref{eq:sigmaInv}, maps null set to null set.
This can be verified by applying Lemma~\ref{lemma:Null1}.

\gap
Theorem~\ref{thm:AkaikeBound} then follows if we can establish the following:

\gap
\noindent
{\bf Claim III}: For almost all $s \in I$, there exists an open set $U_s$ around $\big[\!\begin{smallmatrix} s \\ 0 \end{smallmatrix}\!\big]$
such that $\Theta^{2}(U_s) \subset U_s$.

\gap
\noindent
{\bf Claim IV}: For almost all $p^{(0)} \in {\rm dom}(\Theta)$,
$\Theta^{2k}(p^{(0)}) = {\big[\!\begin{smallmatrix} s \\ 0 \end{smallmatrix}\!\big]}$ with $s \in I$.

\gap
Our proofs of Claims III and IV are based on (essentially part 2 of) the Theorem~\ref{thm:StableManifolds}.

Note that Claim III is local in nature, and follows immediately from Theorem~\ref{thm:StableManifolds} -- also a local result -- for any $s$ in the interior of $I$
excluding $\{\alpha_i, 1-\alpha_i: i=2, \ldots, n-1 \}$. (If we invoke the refined version of Theorem~\ref{thm:StableManifolds}, there is no need to
exclude the singularities. Either way suffices for proving Claim III.)
Claim IV, however, is global in nature.
Its proof combines the refined version of Theorem~\ref{thm:StableManifolds}
with arguments exploiting the diagonal and polynomial properties of $\Theta$.

\gap
\noindent
{\bf Proof of Claim IV.}\footnote{By \eqref{eq:SpectrumPhi},
$\Theta^2$ is a local diffeomorphism at $\big[\!\begin{smallmatrix} s \\ 0 \end{smallmatrix}\!\big]$
for every $s \in (-1,1) \backslash \{\alpha_i, 1-\alpha_i: i=2, \ldots, n-1 \}$.
The interval $I$ defined by \eqref{eq:AttractingInterval}
covers at least 70\% of $(0,1)$: $I \supseteq \big[\frac{1}{2}-\frac{1}{2\sqrt{2}},\frac{1}{2}+\frac{1}{2\sqrt{2}} \big]$; it is easy to check that
$\alpha_{i^\ast}, 1-\alpha_{i^\ast} \in I$, but when $i \neq i^\ast$, $\alpha_i$, $1-\alpha_i$ may or may not fall into $I$.
If we make the assumption that
\bea \label{eq:SingularitiesUnderTheRug}
\alpha_i,\; 1-\alpha_i \in I, \;\; \forall i\neq i^\ast,
\eea
meaning that \emph{all} -- not just one -- the intermediate eigenvalues $\lambda_i$  stay away from $\lambda_1$ or $\lambda_n$,
then part 2 of Theorem~\ref{thm:StableManifolds} can be applied verbatim to $\Theta^2$ at $\big[\!\begin{smallmatrix} s \\ 0 \end{smallmatrix}\!\big]$
for every $s \in (-1,1) \backslash I$ and our argument for proving claim IV will prove the claim under the additional assumption \eqref{eq:SingularitiesUnderTheRug},
and hence a weaker version of Theorem~\ref{thm:AkaikeBound}. Fortunately, Theorem~\ref{thm:StableManifolds} holds without the diffeomorphism assumption -- see the remarks in Appendix~\ref{AppendixB}, so the extra assumption \eqref{eq:SingularitiesUnderTheRug} is not needed.} By \eqref{eq:AkaikeTheorem12}, it suffices to show that the set
$$
\mbox{$\bigcup_{s \in (-1,1)\backslash I}$} \Big\{ p \in {\rm dom}(\Theta): \mbox{$\lim_{k \goto \infty}$} \Theta^{2k}(p) =
\big[\!\begin{smallmatrix} s \\ 0 \end{smallmatrix}\!\big] \Big\}
$$
has measure zero
in $\bR^{n-1}$.
By (the refined version of) Theorem~\ref{thm:StableManifolds} and Lemma~\ref{lemma:Spectrum}, every fixed point
$\big[\!\begin{smallmatrix} s \\ 0 \end{smallmatrix}\!\big]$, $s \in (-1,1)\backslash I$,
of $\Theta^2$
has a center-stable manifold, denoted by $W^{\rm cs}_{\rm loc}(s)$, with \emph{co-dimension at least 1}.
The diagonal property of $T$, and hence of $\Theta^2$, ensures that $W^{\rm cs}_{\rm loc}(s)$
can be chosen to lie on the plane $\{x_i=0: \mu_i(s) > 1\}$, 
which is contained in the hyperplane $\mathcal{P}_\ast := e_{i^\ast}^\perp$. Therefore,
$$
\mbox{$\bigcup_{s \in J}$} W^{\rm cs}_{\rm loc}(s) \subset \mathcal{P}_\ast.
$$
Of course, we also have $\big[\!\begin{smallmatrix} \alpha_i \\ 0 \end{smallmatrix}\!\big]$, $\big[\!\begin{smallmatrix} 1-\alpha_i \\ 0 \end{smallmatrix}\!\big] \in \mathcal{P}_\ast$.
So to complete the proof, it suffices to show that the set of points attracted to the hyperplane $\mathcal{P}_\ast$ by $\Theta^2$ has measure 0, i.e. it suffices to show that
\bea \label{eq:NullSet}
\mbox{$\bigcup_{n\geq 0}$} \Theta^{-2n}(\mathcal{P}_\ast \cap {\rm dom}(\Theta))
\eea
is a null set.

We now argue that $D\Theta^2|_p$ is non-singular for almost all $p \in {\rm dom}(\Theta)$. By the chain rule,
it suffices to show that $D\Theta|_p$ is non-singular for almost all $p \in {\rm dom}(\Theta)$. Note that the entries
of $\Theta(p)$
 are rational functions; in fact $\Theta_i(p)$ is of the form $t_i(p)/v(p)$, where $t_i(p)$ and $v(p)$ are degree
 2 polynomials in $p$ and $v(p)>0$ for $p \in {\rm dom}(\Theta)$. So
 $\det(D\Theta|_p)$ is of the form $w(p)/v(p)^{2(n-1)}$ where $w(p)$ is some polynomial in $p$.
 It is clear that $w(p)$ is not identically zero, as that would violate the invertibility of
$D\Theta|_p$ at many $p$, as shown by \eqref{eq:DTheta_FixedPoint}. It then follows from Lemma~\ref{lemma:Null2}
that $w(p)$ is non-zero almost everywhere, hence
the almost everywhere invertibility of $D\Theta|_p$.

As $\mathcal{P}_\ast \cap {\rm dom}(\Theta)$ is null,
we can then use Lemma~\ref{lemma:Null1} inductively to conclude that
$\Theta^{-2n}(\mathcal{P}_\ast \cap {\rm dom}(\Theta))$ is null for any $n\geq 0$.
So the set \eqref{eq:NullSet}
is a countable union of null sets, hence also a null set.

We have completed the proof of Claim IV, hence also of Theorem~\ref{thm:AkaikeBound}, and Theorem~\ref{thm:Main}(ii) follows. \eop

\gap
\noindent
{\bf Computational examples in 3- and 4-D.} Corresponding to the \textbf{limit probability} $s(x^{(0)})$ for $x^{(0)} \in \bS^{n-1}$, defined by
\eqref{eq:AkaikeTheorem12}, is the \textbf{limit angle} $\theta(x^{(0)}) \in [0,\pi/2]$
defined by $|\hat{x}^{(\infty)}|=[\cos(\theta),0,\ldots,0,\sin(\theta)]^T$, where
$\hat{x}^{(\infty)} := \lim_{k \goto \infty} \GD^{2k}(x^{(0)})/\| \GD^{2k}(x^{(0)}) \|$.
The limit probability and the limit angle are related by
\bea \label{eq:LimitAngles}
\theta = \tan^{-1} \big(a^{-1}\sqrt{s/(1-s)}\big).
\eea\footnote{By Proposition~\ref{prop:Tinvariance}, the map $T$ and (consequently) the limit probability do not depend on ${\rm cond}(A)=a^{-1}$, but the limit angle does.}
This is just the inverse of the bijection between $\theta$ and $s$ \eqref{eq:theta_s} already seen in the 2-D case.
Unlike the 2-D case, an initial vector $x^{(0)}$ uniformly sampled on the unit sphere $\bS^{n-1}$ will not result in
a limit angle uniformly sampled on the interval $[0,\pi/2]$.
In the proof, we establish that $s(x^{(0)}) \in I$ for almost all $x^{(0)} \in \bS^{n-1}$. Therefore,
\bea \label{eq:Interval_J}
\theta(x^{(0)}) \in J := \{\tan^{-1} \big(a^{-1}\sqrt{s/(1-s)}\big) : s \in I  \}, \;\;\; \mbox{$I$ defined in \eqref{eq:AttractingInterval},}
\eea
for almost all $x^{(0)} \in \bS^{n-1}$. The limit angle $\theta(x^{(0)})$ determines the ROC with
the initial vector $x^{(0)}$ according to \eqref{eq:ROC}; this relation is shown by the (light blue) curves in Figure~\ref{fig:Figure4}. (These are the same curves in
the left panel of Figure~\ref{fig:WorstvsAverage}.)

We consider various choices of $A$ with $n=3$, and for each we estimate the
probability density function of $\theta(U)$ with $U\sim {\rm Uniform}(\mathbb{S}^{n-1})$.
As we see from Figure~\ref{fig:Figure4}, the computations suggest
that when the intermediate eigenvalue $\lambda_2$ equals $(\lambda_1+\lambda_3)/2$, the distribution of the limit angles peaks at
$\tan^{-1}(a^{-1})$, the angle that corresponds to the slowest ROC $(1-a)/(1+a)$.
The mere presence of an intermediate eigenvalue, as long as it is not too close to $\lambda_1$ or $\lambda_3$,
concentrates the limit angle $\theta$ near $\tan^{-1}(a^{-1})$. Moreover, the effect gets more prominent
when $a^{-1}={\rm cond}(A)$ is large.
The horizontal lines in Figure~\ref{fig:Figure4} correspond to Akaike's lower bound of ROC in \eqref{eq:LowerBound};
the computations illustrates that the bound is tight.

Finally, we compute the sample minimum, mean and maximum of the ROC based on the four sets of eigenvalues $\lambda$ in Figure~\ref{fig:Figure4}
and two additional $\lambda$'s in 4-D, using in each case a sample of ${\rm Uniform}(\mathbb{S}^{n-1})$ of size $10^7$ as the initial vector.
As we can already tell from the densities in Figure~\ref{fig:Figure4}, the sample mean is significantly closer to the worst case ROC than to the
essential best case ROC.

\begin{center}
\begin{tabular}{|c|c|c|c|c|c|}
  \hline
 $\lambda$  & Essential inf \eqref{thm:AkaikeBound} &  Sample min &  Sample mean &  Sample max & Worst ROC \eqref{eq:AverageWorst} \\
   \hline
$[1,2/3,1/3]$ & 0.37796 &  0.37814 & 0.48618 & 0.50000 & 0.50000 \\
$[1,1/2,1/3]$ & 0.33333 &  0.33360 & 0.46137 & 0.50000 & 0.50000 \\
$[1,5/9,1/9]$ & 0.68599 &  0.68614 & 0.79036 & 0.80000 & 0.80000 \\
$[1,7/9,1/9]$ & 0.63246 &  0.63272 & 0.76854 & 0.80000 & 0.80000 \\
$[1, 4/5, 2/3, 1/3]$ & 0.37796 & 0.37810 & 0.48891 & 0.50000 & 0.50000 \\
$[1, 4/5, 2/5, 1/3]$ & 0.35044 & 0.35046 & 0.47128 & 0.50000 & 0.50000 \\
  \hline
\end{tabular}
\end{center}

\section{Conclusion and Open Questions}

In Section~\ref{sec:POP}, we observe that a stylized version of the extensively studied phase retrieval problems
(see the review article \cite{FannjiangStrohmer:PhaseRetrieval}) is a \emph{well-conditioned} degree 4 POP that can be solved effectively using a \emph{tailored} exact line search GD
algorithm which is only 50\% more expensive per iteration than its constant step size counterpart.
The use of exact line search promises a rate of convergence only matched by the optimally tuned (constant) step size \eqref{eq:OptimalConstant},
which can almost never be achieved in practice.
The theoretical results in \cite{MR3332993} suggests that GD methods are effective
also for two other degree 4 POPs arising from data sciences, namely, the matrix completion problem and the blind source separation problem.
Exact line search algorithms similar to the one proposed in Section~\ref{sec:POP} for phase retrieval
 can be derived for these other problems also.
In our ongoing work, we also learn that the proposed Algorithm~\ref{alg:ExactLineSearchPhaseRetrieval} in Section~\ref{sec:POP} can be adapted to the more realistic setting of phase retrieval imaging with coded diffraction patterns.


Motivated by these applications of exact line search GD methods, we revisit the classical rate of convergence problem of the
method.
In the case of quadratic objectives, the worst case ROC analysis amounts to the Kantorovich inequality, while
the analysis for the average and essential best case ROC requires the analysis of a specific dynamical system.
The latter was spearheaded in the elegant paper of Akaike \cite{Akaike-1959} from more than 65 years ago.
We fill in a gap in this work, using a dynamical system result not available at the time \cite{Akaike-1959} was written.

The computational results in Appendix~\ref{AppendixC} suggest that our main results
Theorem~\ref{thm:Main} and \ref{thm:AkaikeBound} should extend to objectives beyond strongly convex quadratics.
The relatively recent article \cite{deKlerk2017} shows how to extend the
worst case ROC \eqref{eq:Worst} for strongly convex quadratics to general strongly convex functions:
\begin{theorem}\cite[Theorem 1.2]{deKlerk2017}
If $f \in C^2(\bR^n, \bR)$ satisfies $L\cdot I \succeq \nabla^2 f(x) \succeq \mu \cdot I$
for all $x \in \bR^n$, $x^\ast$ a global minimizer of $f$ on $\bR^n$, and $f^\ast = f(x^\ast)$.
Each iteration of the GD method with exact line search satisfies
\bea
f\big( x^{(k+1)} \big) - f^\ast \leq \Big( \frac{L-\mu}{L+\mu} \Big)^{2} \Big( f\big( x^{(k)} \big) - f^\ast \Big), \;\; k=0,1,\ldots.
\eea
\end{theorem}
The result is proved for the slightly bigger family of strongly convex functions $\mathcal{F}_{\mu,L}(\bR^n)$ that requires only $C^1$ regularity; see \cite[Definition 1.1]{deKlerk2017}
or \cite[Chapter 7]{Beck:Book1} for the definition of $\mathcal{F}_{\mu,L}(\bR^n)$.

The above result was anticipated by most optimization experts; compare it with, for example, \cite[Theorem 3.4]{NoceWrig06}.
But the rigorous proof in \cite{deKlerk2017} requires a tricky semi-definite programming formulation; the SDP is designed to
find the initial point that yields the \emph{slowest} ROC.

Extending the main result of this article to objective functions in $\mathcal{F}_{\mu,L}(\bR^n)$ necessarily involves dynamical system
ideas, and necessarily involves a new idea. A special feature of the quadratic case is that the $\GD$ map \eqref{eq:GD} decouples into
the shrinking factor function $\rho$ in $\eqref{eq:rho}$ and the map $[\GD]$ in \eqref{eq:Tabstract}; recall Proposition~\ref{prop:ROC_T}.
This decoupling does not hold for a general objective function, hence an obvious technical challenge.


\appendix
\section{The non-finite termination property of $\GD$}
The following shows that $[\GD]$ is a well-defined self-map on $\mathcal{M}$, where $\mathcal{M}$ is defined in \eqref{eq:NaturalDomain}.
It also shows that unless the initial vector is an eigenvector of $A$, the exact line search GD method cannot reach the minimizer in a finite number of steps.
Compare it with, for example, the finite termination property of the conjugate gradient method.
\begin{proposition}[{\bf ``Non-finite termination"}] \label{prop:Nonzero}
For any $x\neq 0$,
$\GD(x)=0$ if and only if $x$ is an eigenvector of $A$.
If $\GD(x) \neq 0$, then $\GD^k(x)\neq 0$ for all $k$.
\end{proposition}
\pf
Note that $\GD(x)=0$ means $A x = (\frac{x^T A^3 x}{x^T A^2 x}) x$, which is equivalent to $x$ being an eigenvector of $A$.
This establishes the first statement. To prove the second statement, it suffices to show that $\GD(x) \neq 0$ implies $\GD(\GD(x)) \neq 0$.

Assume that the distinct eigenvalues of $A$ partition $\{1,\ldots,n\}$ into $J_1,\ldots,J_m$, i.e.
$\lambda_p = \lambda_q$ iff $p,q \in J_j$ for some $j$. The eignespaces of $A$ are then
$V_\ell={\rm span}(e_i: i \in J_\ell)= \{x\in \bR^n: x_i=0, i\notin J_\ell \}$, $\ell=1,\ldots,m$.
Assume $\GD(x) \neq 0$; by the first statement it is equivalent to that $x \notin V_j$ for all $j$.
We seek to prove that $\GD(x)$ also cannot lie in any of the eigenspaces $V_j$.

Assume the contrary that $\GD(x) \in V_\ell$. So $\GD(x)_i = 0$, $\forall i \notin J_\ell$.
We claim that $x_i=0$, $\forall i \notin J_\ell$, or $x \in V_\ell$, which causes a contradiction.

We argue that
\bea \label{eq:1}
x_i=0 \;\; \forall i\in J_{\ell'}, \;\; \ell'<\ell,
\eea
the proof for $\ell'>\ell$ is similar.
Fix a $\overline{\ell}$ less than $\ell$. (If there is none, there is nothing to show.)
Assume that we have already proved \eqref{eq:1} for all $\ell'<\overline{\ell}$. (If $\overline{\ell}=1$, this is automatic.)
Then, by \eqref{eq:GD},
$(\sum_{j=1}^n \lambda_j^3 x_j^2) x_i - (\sum_{j=1}^n \lambda_j^2 x_j^2)  \lambda_i x_i = 0$, 
$i \in J_{\overline{\ell}}$.
It follows that
\bea \label{eq:Main}
x_i \Bigg\{ \sum_{j\in \bigcup_{k<\overline{\ell}} J_k } \lambda_j^2 (\lambda_j - \lambda_i) \underbrace{x_j^2}_{=0} +
\sum_{j \in J_{\overline{\ell}} } \lambda_j^2 \underbrace{(\lambda_j - \lambda_i)}_{=0} x_j^2 +
\sum_{j\in \bigcup_{k>\overline{\ell}} J_k} \underbrace{\lambda_j^2 (\lambda_j - \lambda_i)}_{<0} x_j^2 \Bigg\}= 0, \;\; \forall i \in J_{\overline{\ell}}.
\eea
If $x_i \neq 0$ for some $i \in J_{\overline{\ell}}$, then \eqref{eq:Main} implies $x_j=0$ for all $j\notin J_{\overline{\ell}}$,
which means $x \in V_{\overline{\ell}}$, a contraction.
Therefore, $x_i = 0$ for $i \in J_{\overline{\ell}}$. We have proved inductively that
$x_i = 0$ for $i \in J_1 \cup \cdots \cup J_{\ell'}$ for any $\ell'<\ell$, and \eqref{eq:1} is proved.
\eop

\section{Two Auxiliary Measure Theoretic Lemmas.}  \label{AppendixA}
The proof of Theorem~\ref{thm:AkaikeBound} uses the following two lemmas.
\begin{lemma} \label{lemma:Null1}
Let $U \subset \bR^m$ and $V\subset \bR^n$ be open sets, $m\geq n$,
 $f: U \goto V$ be $C^1$ with ${\rm rank} (Df(x)) =n$ for almost all $x$.
Then whenever $A$ is of measure zero, so is $f^{-1}(A)$.
\end{lemma}

\begin{lemma} \label{lemma:Null2}
The zero set of any non-zero polynomial has measure zero.
\end{lemma}
For the proofs of these two lemmas, see \cite{3215996} and \cite{1920302}.

\section{The Theorem of Center and Stable Manifolds}  \label{AppendixB}
\begin{theorem}[Center and Stable Manifolds] \label{thm:StableManifolds}
Let 0 be a fixed point for the $C^r$
local diffeomorphism $f: U \goto \bR^n$ where $U$ is a neighborhood of zero in $\bR^n$ and
$\infty > r \geq 1$. Let $E^{\rm s} \oplus E^{\rm c} \oplus E^{\rm u}$ be the invariant splitting of $\bR^n$ into the generalized
eigenspaces of $Df(0)$ corresponding to eigenvalues of absolute value less
than one, equal to one, and greater than one. To each of the five $Df(0)$ invariant
subspaces $E^{\rm s}$, $E^{\rm s} \oplus E^{\rm c}$, $E^{\rm c}$, $E^{\rm c} \oplus E^{\rm u}$, and $E^{\rm u}$ there is associated a local f invariant $C^r$
embedded disc $W^{\rm s}_{\rm loc}$, $W^{\rm cs}_{\rm loc}$, $W^{\rm c}_{\rm loc}$, $W^{\rm cu}_{\rm loc}$, $W^{\rm u}_{\rm loc}$
tangent to the linear subspace at $0$ and a ball $B$ around zero in a (suitably defined) norm such that:
\begin{enumerate}
  \item $W^{\rm s}_{\rm loc} = \{ x\in B | \mbox{ $f^n(x) \in B$ for all $n\geq 0$  and $d(f^n(x),0)$ tends to zero exponentially} \}$.
        $f: W^{\rm s}_{\rm loc} \goto W^{\rm s}_{\rm loc}$ is a contraction mapping.
  \item $f(W^{\rm cs}_{\rm loc}) \cap B \subset W^{\rm cs}_{\rm loc}$. If $f^n(x) \in B$ for all $n \geq 0$, then $x \in W^{\rm cs}_{\rm loc}$.
  \item $f(W^{\rm c}_{\rm loc}) \cap B \subset W^{\rm c}_{\rm loc}$. If $f^n(x) \in B$ for all $n \in \bZ$, then $x \in W^{\rm c}_{\rm loc}$.
  \item $f(W^{\rm cu}_{\rm loc}) \cap B \subset W^{\rm cu}_{\rm loc}$. If $f^n(x) \in B$ for all $n \leq 0$, then $x \in W^{\rm cu}_{\rm loc}$.
  \item $W^{\rm u}_{\rm loc} = \{ x\in B | \mbox{ $f^n(x) \in B$ for all $n \leq 0$  and $d(f^n(x),0)$ tends to zero exponentially} \}$.
        $f^{-1}: W^{\rm u}_{\rm loc} \goto W^{\rm u}_{\rm loc}$ is a contraction mapping.
\end{enumerate}
\end{theorem}
The assumption that $f$ is invertible in Theorem~\ref{thm:StableManifolds} happens to be unnecessary.
The proof of existence of $W^{\rm cu}_{\rm loc}$ and $W^{\rm u}_{\rm loc}$,
based on \cite[Theorem III.2]{MR869255}, clearly does not rely on the invertibility of $f$. That for
$W^{\rm cs}_{\rm loc}$ and $W^{\rm s}_{\rm loc}$, however, is based on the applying \cite[Theorem III.2]{MR869255} to $f^{-1}$;
and \emph{it is the existence of $W^{\rm cs}_{\rm loc}$ that is needed in our proof of Theorem~\ref{thm:AkaikeBound}}.
Fortunately, by a finer argument outlined in \cite[Exercise III.2, Page 68]{MR869255} and \cite{443156},
the existence of $W^{\rm cs}_{\rm loc}$ can be established without assuming the invertibility of $f$.
Thanks to the refinement, we can proceed with the proof without the extra assumption \eqref{eq:SingularitiesUnderTheRug}.

\section{The 2-D Rosenbrock function} \label{AppendixC}
While Akaike did not bother to explore what happened to the optimum GD method in 2-D,
the 2-D Rosenbrock function $f(x) = 100(x_2-x_1^2)^2 + (1-x_1)^2$, again a degree 4 polynomial, is often used in optimization textbooks to exemplify different
optimization methods. In particular, due to the ill-conditioning near its global minimizer $x^\ast = [1,1]^T$ (${\rm cond}(\nabla^2 f(x^\ast)) \approx 2500$), it is often used to
illustrate the slow convergence of GD methods.

\begin{figure}[ht]
\begin{center}
\begin{tabular}{ccc}
\includegraphics[height=4.5cm]{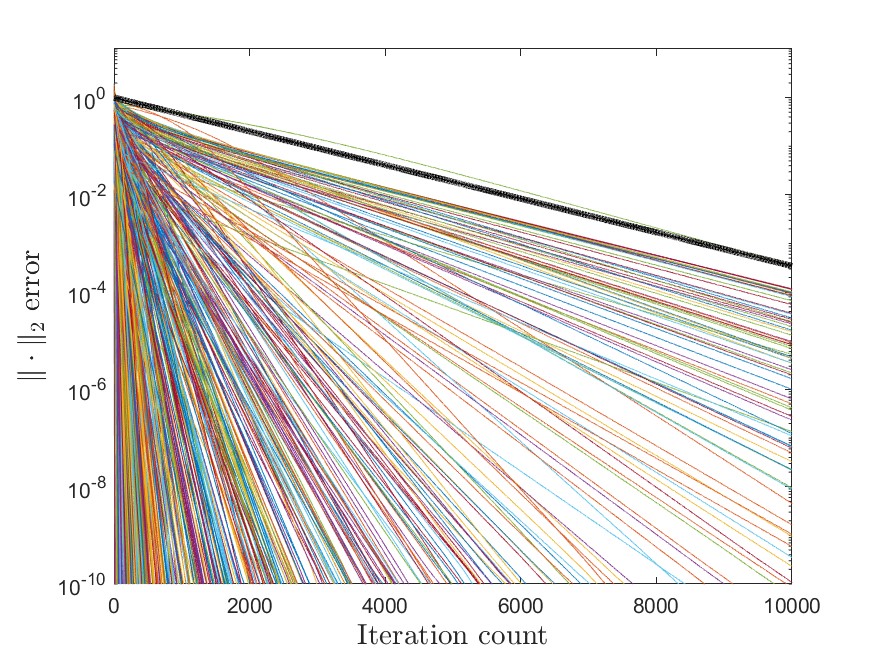} \hspace{-.5cm}
&
\includegraphics[height=4.5cm]{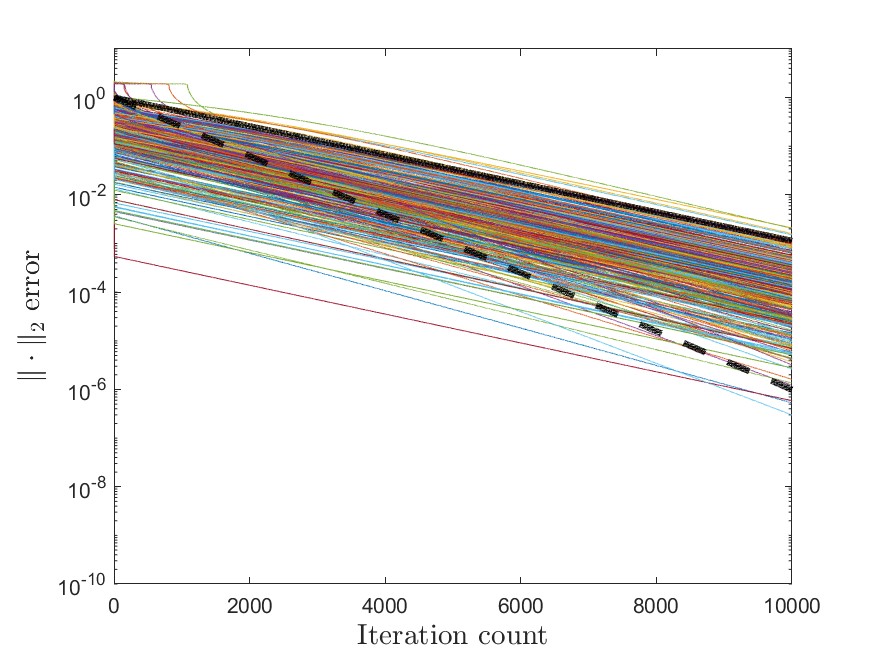} \hspace{-.5cm}
&
\includegraphics[height=4.5cm]{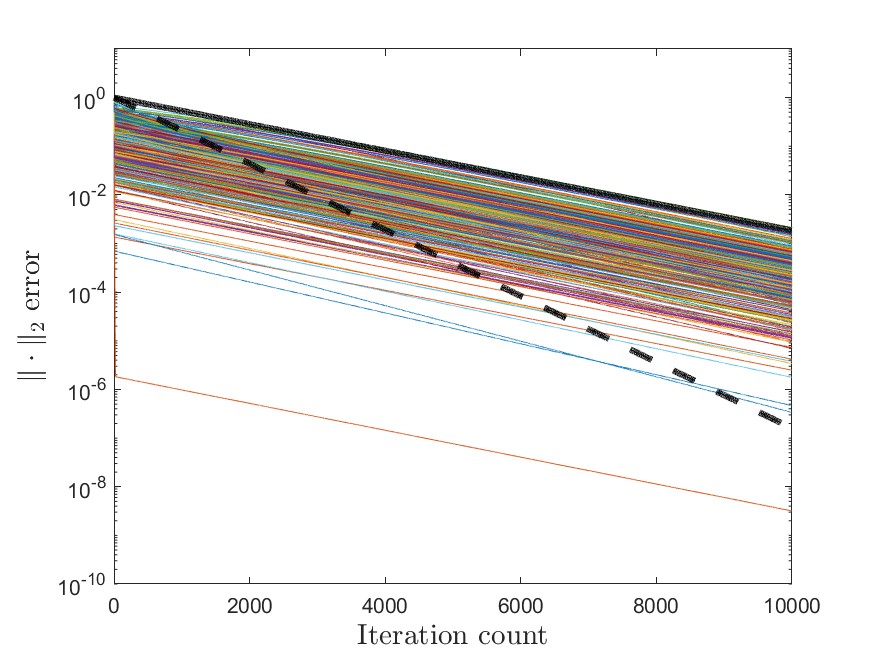} \\
$n=2$ & $n=3$ & $n=4$ \\
\end{tabular}
\end{center}
\vspace{-.6cm}
\caption{ROC of the optimum GD method applied to the Rosenbrock function with $n$ variables with 500 initial guesses uniformly sampled from the unit ball around
the unique minimizer.
The black solid and black dashed lines illustrate the worst case and essential best ROCs, respectively, assuming that the objective were a quadratic with Hessian $\nabla^2 f(x^\ast)$.}
\label{fig:Figure5}
\end{figure}

A student will likely be confused if he applies the exact line search GD method to this objective function.
 Figure~\ref{fig:Figure5} (leftmost panel) shows the ROC of
optimum GD applied to the 2-D Rosenbrock function. The black line illustrates the worst case ROC given by \eqref{eq:AverageWorst}
under the pretense that the degree 4 Rosenbrock function is a quadratic with the constant Hessian $\nabla^2 f(x^\ast)$;
the other lines show the ROC
with 500 different initial guesses
sampled from $x^\ast + z/\|z\|_2$, $z \sim N(0,I_2)$. As predicted by the first part of Theorem~\ref{thm:Main}, the average ROC is much
faster than the worst case ROC shown by the black line not in spite of, but because of, the ill-conditioning of the Hessian.\footnote{To the best of the author's knowledge,
the only textbook that briefly mentions this difference between dimension 2 and above is
\cite[Page 62]{NoceWrig06}. The same reference points to a specialized analysis for the 2-D case in an older optimization textbook, but the latter does not contain a result to the effect of
the first part of Theorem~\ref{thm:Main}.}

The student may be less confused if he tests the method on the higher dimensional Rosenbrock function:
$$
f(x) = \sum_{i=2}^n 100 (x_i - x_{i-1}^2)^2 + (1-x_{i-1})^2.
$$
See the next two panels of the same figure for $n=3$ and $n=4$, which are consistent with what the second part of Theorem~\ref{thm:Main} predicts, namely, ill-conditioning
leads to slow convergence for essentially all initial vectors.



\bibliographystyle{plain}
  {\small
\bibliography{refinement}
  }
\end{document}